\documentclass[12pt]{article}    
\usepackage{amsmath}   
\usepackage{amsfonts}   
\usepackage{amssymb}   \usepackage{times}

    \def\d{\mathop{\hbox{\rm d}}}   
   \def\Diff{\mathop{\hbox{\rm Diff}}}

  \newtheorem{prob}{Problem}    
\newtheorem{prop}{Proposition}[section]    \newtheorem{coro}{Corollary}[section]   
\newtheorem{lem}{Lemma}[section]    \newtheorem{theo}{Theorem}[section]   
\newtheorem{defi}{Definition}[section]    
    
\newtheorem{cri}{Sufficient condition}       
    
\newtheorem{rem}{Remark}[section]   \newtheorem{exa}{Example}[section]    
\newtheorem{hyp}{Hypothesis}[section]

\newfont{\gothic}{eufm10 at 12pt}     
\title{The 1856  Lemma of Cayley Revisited}    
\author{Jin Yue\footnote{E-mail: yue@mathstat.dal.ca}\\      
Department of Mathematics and Statistics\\    
Dalhousie University, Halifax, Nova Scotia \\  Canada B3H 3J5 }     
\newcommand{\D}{\displaystyle}  
\begin{document}    \maketitle   \setlength{\baselineskip}{22pt}      
\begin{abstract} The result of the classical invariant theory (CIT) commonly   
referred to as Lemma of Cayley is reviewed. Its analogue in the   invariant 
 theory of Killing tensors (ITKT) defined in pseudo-Riemannian   spaces of constant 
curvature is formulated and proven. Illustrative examples  are provided.     \end{abstract} 
\medskip  \noindent {\em Subj. Class.:} Differential geometry  \medskip  
   \bigskip  \section{Introduction}   
 In recent years the classical invariant theory (CIT)  of  homogeneous  polynomials has 
reinvented itself once again through new aspects  of the  Lie group theory 
(notably, the generalizations of the moving frames method   due to Fels and 
Olver \cite{FO1, FO2} and Kogan \cite{KO},  see also the relevant  references therein),
 the rise of the modern  computer algebra and new applications  in other areas 
of mathematics  (see Hilbert \cite{Hil} and Olver \cite{PJO}  for a complete review  
and related references). Thus, in  their pioneering 2002 paper  
 McLenaghan {\em et al} \cite{MST1}  successfully planted  the underlying ideas of CIT  
into the fertile field of the  (geometric) study of Killing tensors defined in 
 pseudo-Riemannian manifolds  of constant curvature, which ultimately bore the  
fruit of a new theory  (see also \cite{HMS, SY,  MST2, DHMS, MST3, MST4, MST5, MST6}).   
The resulting  {\em invariant theory of Killing tensors (ITKT)}  shares many of  
the same essential features with the  original CIT.  In light of the fact  that 
``Mathematics is the study of analogies  between  analogies'' \cite{GCR}, 
 we wish to continue developing  ITKT  by establishing more analogies with  CIT.       
  As is well known, the {\em main object} of study in  CIT is a  vector
 space of homogeneous  polynomials  under the action of the  general linear group
 (or its subgroups),  while the {\em main problem}  is that of the determination of  
the functions of the  parameters of the vector  space in question that remain fixed 
under the action of the group.   These functions, called  {\em invariants} 
(Sylvester is credited as the first  to coin  the term), are very useful  in solving  
various classification problems.     In this study the vector spaces of particular 
 importance are the  spaces of {\em binary forms}, or homogeneous polynomials 
 of degree $n$ in two variables,  originally referred to by Cayley  as {\em quantics}. 
 Let ${\cal Q}^n(\mathbb{R}^2)$ denote the vector space  of binary forms of 
 degree $n$ over the reals. Then the dimension $d$ of  the space is given by    
\begin{equation}   d=\dim\,{\cal Q}^n(\mathbb{R}^2) = n + 1.  \label{dQ}   \end{equation} 
 The  general form of an element ${Q}(x,y)$ of the vector space   
 ${\cal Q}^n(\mathbb{R}^2)$ is determined by  the following formula.   
\begin{equation}   {Q}(x,y)=\sum_{i=0}^{n}{n\choose i} a_ix^{n-i}y^i,   
\quad (x,y) \in \mathbb{R}^2. \label{gQ}   \end{equation}   
 Note the arbitrary constants  $a_0, \ldots, a_{n}$ represent the parameter   
space $\Sigma \simeq \mathbb{R}^{n+1}$ corresponding to ${\cal Q}^n(\mathbb{R}^2)$.    
The  special linear group $SL(2,\mathbb{R})$ (for example)  acts on the space  
${\cal Q}^n(\mathbb{R}^2)$ by linear substitutions,   which yield the corresponding  
transformation rules    
\begin{equation}  \begin{array}{c}   
\tilde{a}_0 = \tilde{a}_0(a_0,\ldots, a_n, \alpha, \beta, \gamma, \delta),  \\ 
 \tilde{a}_1 = \tilde{a}_1(a_0,\ldots, a_n, \alpha, \beta, \gamma, \delta),  \\ 
 \vdots \\  \tilde{a}_n = \tilde{a}_n(a_0,\ldots, a_n, \alpha, \beta, \gamma, \delta), 
 \\   \end{array} \label{TLP}  
\end{equation}  
 where  $\alpha, \beta,\gamma, \delta \in \mathbb{R},$  
 $\alpha\delta - \beta\gamma = 1$ 
are local coordinates that parametrize  the group.  Note $\dim\, SL(2,\mathbb{R}) = 3.$ 
The formulas (\ref{TLP})  can be derived  explicitly \cite{PJO}.   
The problem is now reduced to  finding all of the  {\em invariants} of the 
$SL(2,\mathbb{R})$ action  on the space $\Sigma$,  or the functions of 
$a_0,\ldots, a_{n}$ that remain  unchanged under the  transformations (\ref{TLP}):  
 \begin{equation}   
{\cal I} = F(\tilde{a}_0, \ldots, \tilde{a}_{n}) = F(a_0, \ldots, a_{n}). 
 \label{INV}   \end{equation}  
 Note that in the case of $SL(2,\mathbb{R})$ acting on the vector space  
the invariants  appear to be of weight zero due to the condition 
 $\alpha\delta - \beta\gamma = 1$. 
 In order to describe the space  of all $SL(2,\mathbb{R})$-invariants of the vector space 
  ${\cal Q}^n(\mathbb{R}^2)$ one has to  determine a set of the   
{\em fundamental invariants}, with the property that all other invariants  
 are (analytic) functions of the fundamental invariants.   
The number of fundamental  invariants can be determined by using the result 
 of the Fundamental Theorem on  Invariants of  a regular Lie group action \cite{PJO}:  
 \begin{theo}   
Let $G$ be a  Lie group acting regularly on an $m$-dimensional  manifold
 $X$ with $s$-dimensional  orbits.  Then, in a neighborhood $N$  
of each point $x_0 \in X$, there exist $m-s$  functionally independent 
 $G$-invariants $\Delta_1, \ldots,$ $ \Delta_{m-s}$.   
Any other $G$-invariant $\cal I$ defined near $x_0$ can be locally  uniquely expressed 
 as an analytic function of the fundamental invariants   through
 ${\cal I} = F(\Delta_1,$  $ \ldots,$ $ \Delta_{m-s})$.  \label{FT}   
\end{theo}   
The following proposition \cite{PJO}   provides a mechanism for determining
 the dimension of the orbits of  a regular  Lie group action.    
\begin{prop}   Let a Lie group $G$ act on $X$, {\gothic g} is the Lie algebra of 
 $\,G$ and let $x\in X$.   The vector space
 ${ S}|_x = \mbox{\em Span}\{{\bf V}_i(x)|\,   {\bf V}_i \in \mbox{\gothic g}\}$ 
 spanned by all vector fields determined   by the infinitesimal generators at $x$ 
coincides with the tangent space    to the orbit  ${\cal O}_x$ of $G$ 
that passes through $x$,   so $S|_x = T{\cal  O}_x|_x$.
In particular, the dimension of ${\cal O}_x$   equals the  dimension of 
$S|_x$.  \label{Prop1}   
\end{prop}   
One way to determine the fundamental invariants is to use the  infinitesimal 
 generators of the Lie algebra of the  group, by which  we mean their counterparts 
in the parameter space $\Sigma$  satisfying the same commutator relations as the  
generators  defined in the original space. Thus,  a function $F(a_0,\ldots, a_{n})$  
 is an invariant iff it is annihilated by the generators of the Lie algebra  
 defined in the parameter space $\Sigma$. Accordingly, 
 the problem of the  determination of a set of the fundamental invariants 
reduces to solving the  corresponding system of PDEs defined by  the generators. 
This is a short description  of  Sophus Lie's  
{\em method of the infinitesimal generators}, which can be used  to  
compute the invariants. Another powerful method, about  which we shall not dwell   
in this paper, is   \'{E}lie Cartan's   {\em method of moving frames}, 
 which has been recently brought back  to light  \cite{FO1,FO2,PJO, PJO1, BO, BMS, KO, 
SY, DHMS}.    
Arthur Cayley's main contributions to the development of CIT  appeared during  
the period 1854-1878 in his famous ``ten memoirs on quantics''. 
 Having introduced  the notion of an abstract group, he was the first  to recognize 
that the action of a  Lie group on a vector space  can be investigated by studying 
its  `` infinitesimal action'',  that is the corresponding Lie algebra.
 In spite of   the fact that  Cayley thought of this as of something pertinent 
only to   the general linear group and its subgroups, his results  in this area may  
be  considered as a precursor to  Sophus Lie's theory of abstract Lie groups  
that was  developed later in the 19th century.  More specifically,  
in his ``second memoirs on  quantics'' \cite{AC} Arthur Cayley  considers
 (in  modern mathematical language)  the problem of  the determination of 
the action of the Lie group $SL(2,\mathbb{R})$  
 on the vector space ${\cal Q}^n(\mathbb{R}^2)$ in conjunction with  the problem of  
 computing the invariants.  The main result is the subject of the following lemma 
  (see Cayley \cite{AC} and Olver \cite{PJO}, p.213).     
\begin{lem}[Cayley] \label{cayley}    
The action of $SL(2, \mathbb{R})$ on the space  ${\cal Q}^{n}(\mathbb{R}^2)$ of   
 binary homogeneous polynomials of degree $n$  defined by (\ref{gQ}) 
has the following   infinitesimal generators  
in the corresponding parameter space $\Sigma$:   
 \begin{equation} \begin{array}{rcl}   
{\bf V}^{-}   &=&   
na_1\partial _{a_0}+ (n-1) a_2\partial_{a_1}+ \cdots  
 + 2a_{n-1}\partial _{a_{n-2}}  +a_{n}\partial _{a_{n-1}},  \\ [0.3cm] 
 {\bf V}^{0}  &=&   
- n a_0\partial _{a_0} +(2-n)a_1\partial_{a_1}+\cdots  
 + (n-2)a_{n-1}\partial _{a_{n-1}} \\ [0.3cm]  & & + n a_{n}\partial _{a_{n}}, \\ [0.3cm] 
  {\bf V}^{+} &=& a_0\partial _{a_1}+2a_1\partial_{a_2} 
 +\cdots  + (n-1)a_{n-2}\partial _{a_{n-1}}  + n a_{n-1}\partial _{a_{n}},     
\end{array} \label{gen}  \end{equation}  
 where  $\partial _{a_i}=\frac{\partial } {\partial a_i},$  $i=0, \ldots,n$.   
\end{lem}    
Observe that the vector fields (\ref{gen}) enjoy the following  commutator relations   
\begin{equation}   
[{\bf V}^-, {\bf V}^0] = -2{\bf V}^-,\quad   
[{\bf V}^+,  {\bf V}^0] = 2{\bf V}^+, \quad   
[{\bf V}^-, {\bf V}^+] = {\bf V}^0,  \label{com}   
\end{equation}   which confirm that the generators (\ref{gen}) represent the action  
 of $SL(2, \mathbb{R})$ in the parameter space $\Sigma$.   
  In view of the above, solving the problem of the determination of  the 
 $SL(2,\mathbb{R})$-invariants of the vector space  ${\cal Q}^n(\mathbb{R}^2)$
  for a specific $n$ amounts now  to solving the corresponding system of linear 
 PDEs  determined by the generators (\ref{gen}):    
\begin{equation} \begin{array}{lll}   
{\bf V}^- (F) &=& 0, \\ [0.3cm]   {\bf V}^0 (F) &=& 0, \\ [0.3cm]   
{\bf V}^+ (F) &=& 0. \label{gPDE}   \end{array} \end{equation}    
for a (analytic) function $F$ defined in the parameter space $\Sigma$.  
 We note that according to Proposition \ref{Prop1} the dimension of  the orbits of 
 the $SL(2,\mathbb{R})$ action on ${\cal Q}^n(\mathbb{R}^2)$  can (locally) be
  determined by the number of linearly independent vector  fields (\ref{gen}).  
 Accordingly, by Theorem \ref{FT}  the number of fundamental  
$SL(2,\mathbb{R})$-invariants is  $n+1 - s$,  where $s \le 3$ is the dimension   
of the orbits.  Therefore for each particular $n$ the general solution to the system  
 (\ref{gPDE}) will take the form   
\begin{equation}   
{\cal I} = F(\Delta_1, \Delta_2, \ldots, \Delta_{\ell}), \label{gS}   \end{equation}   
where $\Delta_i = \Delta_i(a_0, \ldots, a_{n}),$ $i=1, \ldots, \ell$,   
$\ell = n+1-s$  are the fundamental $SL(2,\mathbb{R})$-invariants.    
 To illustrate the procedure, let us recall the  following  
well-known example  \cite{PJO}.  
\begin{exa} {\rm  Consider the  vector space  ${\cal Q}^2(\mathbb{R}^2)$. 
The   elements of  ${\cal Q}^2(\mathbb{R}^2)$ enjoy the following general form. 
\begin{equation} \label{2Q}   Q(x,y)= a_0x^2 + 2a_1xy  + a_2 y^2.   \end{equation}   
The (local) action of $SL(2,\mathbb{R})$ in the parameter space  
 $\Sigma \simeq \mathbb{R}^3$ generated by the parameters $a_0, a_1$ and   
$a_2$ is represented by the  vector fields  
 \begin{equation} \begin{array}{rcl}   
{\bf V}^{-}   &=&  a_0\partial_{a_1} + 2a_1\partial_{a_2},  \\ [0.3cm]  
{\bf V}^{0}   &=&  2a_0\partial_{a_0} - 2a_2\partial_{a_2}, \\ [0.3cm]   
{\bf V}^{+}   &=&  2a_1\partial_{a_0} + a_2\partial_{a_1}.     
\end{array} \label{2dim}  \end{equation}   
obtained via the standard technique of exponentiation. 
We immediately observe that only   two vector fields (\ref{2dim})  
are linearly independent,  therefore in view of Theorem \ref{FT} and 
 Proposition \ref{Prop1}  there is (almost everywhere) $3-2=1$ fundamental 
 $SL(2,\mathbb{R})$-invariant  of the vector space ${\cal Q}^2(\mathbb{R}^2)$. 
 Indeed, solving the system  of PDEs (\ref{gPDE}) for the vector fields (\ref{2dim})  
 yields the solution:  $$ {\cal I} = F(\Delta_1),$$ where $\Delta_1 = a_0a_2 - a_1^2.$  
 The group acts with  orbits of two types:  $a_0 = a_1 = a_2 = 0$, 
which  is an orbit of dimension $0$ and the level sets  of
 $\Delta_1$ (i.e., $\Delta_1 =0$  and $\Delta_1 \not= 0$), 
 both of which are orbits of dimension $2$.    }  \end{exa}   
 Now let us turn our attention to ITKT.  Here the underlying space is  
 a pseudo-Riemannian manifold $(M, {\bf g})$ of constant curvature. 
 The   vector spaces  in question are the vector spaces of Killing tensors. 
 Our notations  are compatible  with those introduced in \cite{MST2}. 
 Thus, for a fixed $n \ge 1,$  ${\cal K}^n(M)$ denotes the vector space 
 of Killing tensors of valence $n$  defined on $(M, {\bf g})$. 
 The group acting on ${\cal K}^n(M)$ is the isometry  group $I(M)$ of  
$(M, {\bf g})$. \begin{rem}{\rm Here and below $I(M)$ denotes the
 {\em continuous} Lie group of isometries of $M$. We do not take into consideration 
discrete isometries. } 
\end{rem}  
A  comprehensive review of ITKT is the subject of  Section 2. 
 Now, let us formulate an analogue of the problem solved by Cayley \cite{AC}.  
 Since Cayley's problem concerns binary forms it will be natural to  investigate  
in this respect the Killing tensors of arbitrary valence  defined in pseudo-Riemannian 
 manifolds of dimension two,  for example, the Minkowski plane $\mathbb{R}_1^2$,
 More infromation about the Minkowski geometry can be found in Thompson \cite{AT}.
  Accordingly, the vector  spaces that we shall study in what follows are  
${\cal K}^n(\mathbb{R}^2_1)$,  $n \ge 1$. Table \ref{Tab1} presents a comparison  
of the ``ingredients''  and information that can be used to solve  the  
two sister-problems.      
\begin{table}[ht]   \begin{center}   \begin{tabular}{|c|c|c|c|c|} \hline   
{\rule[-7mm]{0pt}{16mm}{\small Theory}} & \begin{tabular}{c}  
 {\small Vector} \\[0.1cm] {\small space}  \end{tabular} &  {\small Group} & 
\begin{tabular}{c}  {\small Dimension of}  \\[0.1cm]  {\small the space} \end{tabular}
  & \begin{tabular}{c} {\small Dimension}  \\[0.1cm]  {\small of the orbits} \end{tabular}
 \\[0.2cm]\hline     \begin{tabular}{c} { }\\[0.1cm]    CIT \\[0.1cm] {}   
\end{tabular} &  ${\cal Q}^n(\mathbb{R}^2)$ &  $SL(2,\mathbb{R})$ 
 & $n+1$  &  $ \le 3 $    \\  \hline   \begin{tabular}{c}{ }
\\[0.1cm] ITKT \\[0.1cm] { }   \end{tabular} & ${\cal K}^n(\mathbb{R}^2_1)$ & 
 $I(\mathbb{R}_1^2)$ &   $\displaystyle \frac{1}{2}(n+1)(n+2) $ &    $ \le 3 $ 
 \\ [0.2cm]  \hline  \end{tabular}    \caption{\small 
The settings for the corresponding problems in CIT and ITKT.}    \label{Tab1}    
\end{center}    \end{table}   
 Having made these observations, we are now in the position  
 to formulate  the ITKT version of the problem considered by Cayley  in \cite{AC}.    
\begin{prob}  Consider the action of the isometry Lie group $I(\mathbb{R}^2_1)$   
 on the vector space ${\cal K}^n(\mathbb{R}_1^2)$. 
 Determine  a representation of  the corresponding Lie algebra  
$i(\mathbb{R}_1^2)$ on the parameter space $\Sigma$  of  ${\cal K}^n(\mathbb{R}_1^2)$. 
 \label{Prob}    \end{prob}    Clearly, the solution to this problem will mimic
 the result of  Lemma  \ref{cayley}, namely one will have to determine  a basis of 
the  Lie algebra  defined on the parameter space $\Sigma$ of  
${\cal K}^n(\mathbb{R}^2_1)$,  which is isomorphic  to the  Lie algebra
 $i(\mathbb{R}_1^2)$. Having the generators  of such a Lie algebra 
 will allow one to compute the   $I(\mathbb{R}_1^2)$-invariants of  
 ${\cal K}^n(\mathbb{R}_1^2)$ by solving the corresponding system of PDEs  in the   
 spirit of the corresponding problem of CIT described above.       
To solve the problem we need to establish first  the requisite language of ITKT.  
 This is the subject of   the considerations that follow in Section 2.          

\section{Invariant theory of Killing tensors (ITKT)}      
Perhaps the most efficient  way to begin describing a mathematical theory  
is by placing it among  other  mathematical theories. 
 Recall that in the 19th century the post-``Theorema Egregium  of Gauss"  
differential geometry branched off into two  directions. Thus,   B. Riemann \cite{Ri} 
generalized the theory of surfaces of C. F. Gauss,   from two to  several dimensions, 
which ultimately led  to the emergence  of the  new geometric objects  known
 now as {\em (pseudo-) Riemannian manifolds},   and  more broadly, 
 today's {\em differential geometry}.   The other school of thought was based on
 F. Klein's idea that every geometry   could be interpreted as a theory of invariants 
with respect to a specific   {\em transformation group}. Thus,  
according to F. Klein \cite{FK72, FK93},   the main objective of any branch of 
geometry can be described  as follows:  ``Given a manifold and a group of 
transformations of the manifold,   to study the manifold configurations with respect 
to those features   that are not altered by the transformations of the group'' 
(\cite{FK93}, p.67).   One of the most fundamental contributions of \'{E}. Cartan, 
 in  particular,  with his {\em theory of moving frames} \cite{Car},  
is the fusion of these  two directions into a single theory. The comprehensive 
monograph by Sharpe  \cite{Sha} unveils all of the beauty  of Cartan's theory that
 subsumed the ideas  of both Riemann and Klein (see also, for example, 
Arvanitoyeorgos \cite{Arv}). The following diagram (see \cite{Sha}, p.ix) describes
 the  relationship among the different approaches to geometry  mentioned above.  
\begin{equation}  \begin{array}{ccc}  \mbox{Euclidean Geometry} &  
\stackrel{\mbox{\small generalization}}{\longrightarrow} &  \mbox{Klein Geometries} 
\\[0.5cm]  \downarrow\mbox{\small generalization} &  & 
\mbox{\small generalization}\downarrow \\[0.5cm]  \mbox{Riemannian Geometry} & 
 \stackrel{\mbox{\small generalization}}{\longrightarrow} &  \mbox{Cartan Geometries} 
 \end{array}  \label{Dia}  \end{equation}     
ITKT  can be placed into the
 theory of Cartan linking the developments  of  Riemann and Klein. We now shall 
present the evidence  to justify this claim.   
Indeed, let $(M, {\bf g})$ be a pseudo-Riemannian manifold of  constant curvature, 
 assume also that $\dim\, M = m$.      
\begin{defi}  
A {\em Killing tensor ${\bf K}$  of valence $n$ defined in  $(M,{\bf g})$}
 is a symmetric $(n,0)$-type tensor satisfying  the  Killing tensor equation   
\begin{equation}   [{\bf K}, {\bf g}] = 0, \label{KE}  \end{equation}   
where $[$ , $]$ denotes the Schouten bracket \cite{Sch}.  
 When $n=1$, ${\bf K}$ is said to be a {\em Killing vector (infinitesimal isometry)}  
 and the equation (\ref{KE}) reads $${\mathcal L}_{{\bf K}}{\bf g} = 0,$$  
 where $\mathcal L$ denotes the  Lie derivative operator. \label{DKT}   
\end{defi}     
\begin{rem}  
{\rm  Throughout this paper, unless otherwise specified, 
  $[$ , $]$ denotes the Schouten bracket, which is a generalization of  
the usual Lie  bracket of vector fields.  }  
\end{rem}    
The set of all Killing vectors of $(M,{\bf g})$,   denoted by $i(M)$, 
is a Lie algebra of the corresponding Lie group   of isometries $I(M)$,
 which is also a Lie subalgebra of the space  ${\cal X}(M)$  of all vector fields 
defined on $M$.   As is well-known,  $d = \dim i(M) = \frac{1}{2}m(m+1)$ iff 
the space  $(M, {\bf g})$ is of constant curvature.    It follows immediately 
from (\ref{KE}) that   Killing tensors of  the same valence  $n$ constitute a 
vector space ${\cal K}^n(M)$.  Moreover,  the following properties hold true:   
\begin{eqnarray}   
\mbox{$[$ , $]$:} \,\,\,  
{\cal K}^n(M)\oplus {\cal K}^{\ell} (M)  \rightarrow {\cal K}^{n+\ell-1} (M),  
\\[0.5cm]  [{\bf K}^n, {\bf K}^{\ell}]  = - [{\bf K}^{\ell}, {\bf K}^n] \quad  
\mbox{(skew-symmetry)},  \\[0.5cm] 
 [[{\bf K}^n, {\bf K}^{\ell}], {\bf K}^r] + \mbox{(cycle)} = 0  
 \quad \mbox{(Jacobi identity)}, \label{JI}   
\end{eqnarray}    
where ${\bf K}^n \in {\cal K}^n(M),$  
${\bf K}^{\ell} \in {\cal K}^{\ell}(M),$ ${\bf K}^r \in {\cal K}^r(M).$ 
  Therefore one can consider a graded Lie algebra of Killing tensors  
 defined on $(M, {\bf g})$ with respect to the Schouten bracket  $[$ , $]$:   
\begin{equation}   
{\cal K}_{alg} =  {\cal K}^0(M)\oplus {\cal K}^1(M)\oplus  {\cal K}^2(M)\oplus  
\cdots \oplus {\cal K}^n(M) \oplus \cdots,  
 \end{equation}   
where ${\cal K}^0(M) = {\mathbb R}$,   
 ${\cal K}^1(M) = i(M)$ and $n  = 0, 1, 2, \ldots,$  
 denotes the valence of the Killing tensors belonging to  the corresponding  
space ${\cal K}^n(M)$. These remarkable geometrical  objects have  been  actively 
studied for a long time  by mathematicians and physicists alike. 
 Apart from possessing beautiful  mathematical properties, Killing tensors and  
conformal Killing tensors  naturally arise in many problems of classical mechanics,   
general relativity, field theory and other areas.  More information can be found, 
 for example, in the following references:  Delong \cite{De},  Dolan {\em et al} 
\cite{Do}, Benenti \cite{Be},   Bruce {\em et al} \cite{BMS},  Bolsinov and Matveev
 \cite{BM},   Crampin \cite{Cr},  Eisenhart \cite{Ei1, Ei2},  Fushchich and Nikitin 
\cite{FN},  Kalnins \cite{Kal, Ka},  Kalnins and Miller \cite{KM1, KM2}, 
Miller \cite{Mi},   Mokhov and Ferapontov \cite{MF},  Takeuchi \cite{Ta}, 
Thompson \cite{Th},  as well as many others  (more references related to
the study of Killing tensors  of valence two  can be found in the review \cite{Be}).  
    To illustrate how Killing tensors appear naturally in the problems of  
 classical mechanics, let us consider the following example.  
 \begin{exa}  {\rm   Let $({\bf X}_H, {\bf P}_0, H)$ be a Hamiltonian system defined  
on  $(M, {\bf g})$ by a natural Hamiltonian $H$ of the form   
\begin{equation}   H ({\bf q}, {\bf p}) = \frac{1}{2}g^{ij}p_ip_j + V({\bf q}), 
 \quad i,j  = 1,\ldots, m, \label{H}  
 \end{equation}  
 where $g^{ij}$ are the contravariant components of the corresponding  
 metric tensor $\bf g$,  $({\bf q},{\bf p}) \in T^*M$ are the canonical  
 position-momenta coordinates and the Hamiltonian vector field   ${\bf X}_H$ 
is given by   \begin{equation}   {\bf X}_H = [{\bf P}_0, H] \label{X}   
\end{equation}   
with respect to the canonical Poisson bi-vector   
 ${\bf P}_0 = \sum_{i=1}^m\partial/\partial q^i\wedge \partial/\partial p_i.$  
 Assume also that the Hamiltonian system defined by (\ref{H}) 
 admits a first integral  of motion $F$ which is a polynomial function  
of degree $n$ in the momenta:   
\begin{equation}   
F({\bf q}, {\bf p})=  K^{i_1i_2\ldots i_n}({\bf q})p_{i_1}p_{i_2}\ldots p_{i_n}  
+ U({\bf q}),   \end{equation}  
 where $ 1\le i_1, \ldots, i_n \le m.$  Since the functions 
 $H$ and $F$ are in involution,   the vanishing of the  Poisson bracket defined by  
${\bf P}_0$: $$\{H, F\}_0 = {\bf P}_0\d H\d F = 
[[{\bf P}_0, H], F] = 0$$  yields  \begin{equation}  [{\bf K}, {\bf g}] =0, 
  \quad \mbox{(Killing tensor equation)} \label{KT1}  
 \end{equation}  and   \begin{equation}  
 K^{i_1 i_2\ldots i_n}\frac{\partial V}{\partial q^{i_1}}p_{i_2}\ldots p_{i_n} 
  = g^{ij}\frac{\partial U}{\partial q^i}p_j,   
\quad \mbox{(compatibility condition)}, \label{CC}   
\end{equation}   
where the symmetric $(n,0)$-tensor $\bf K$ has the components  
 $K^{i_1 i_2\ldots i_n}$ and  $1 \le $  $i,$ $j,$ $i_1,$ $\ldots,$ $ i_n$  
 $ \le m$. Clearly, in view of Definition \ref{DKT} the equation   (\ref{KT1})
 confirms that $\bf K$ is a Killing tensor.  Furthermore,  in the case $n=2$ 
(see Benenti \cite{Be})  the compatibility condition (\ref{CC})  reduces to  
${\bf K}\d V = {\bf g} \d U$ or $\d (\hat{\bf K}\d V) = 0$,   
where the $(1,1)$-tensor $\hat{\bf K}$ is  given by  
 $\hat{\bf K} = {\bf K}{\bf g}^{-1}$.   } \label{E1}  \end{exa}   
 Example \ref{E1} ellucidates the appearence of Killing tensors  in the problems  of 
 the integrability theory of Hamiltonian systems.  Notably, the  geometric properties
of Killing tensors of valence two  have been routinely  employed for a long time to solve 
   the problems arising in the  theory of   orthogonal separation of variables  
 (see, for example,  \cite{Bo, Dar, Ei1, Ei2, De,  Mi, Kal, Ka,  Be,  BMS, MST1, MST3, HMS}
  and the relevant references therein).    
Recall that the standard approach to the study of Killing tensors  
 defined in pseudo-Riemannian spaces of constant curvature  
is based on the property  that the Killing tensors defined  in these spaces are
 sums of symmetrized tensor  products   of Killing vectors (see, for example, 
\cite{De, Kal, Th}).        
\begin{exa}  {\rm  Consider the set 
${\cal K}^2(\mathbb{R}^2_1)$ of all Killing tensors   of valence two defined in 
$\mathbb{R}_1^2$ (Minkowski plane).    Recall  that the Lie algebra $i(\mathbb{R}_1^2)$
  of Killing vectors  (infinitesimal isometries) admits  the basis given by the
 following Killing vectors:   
\begin{equation}  {\bf T} = \partial _t,  \quad  {\bf X} = \partial_ x, 
 \quad  {\bf H} = x \partial _t + t \partial _x  \label{TXH}  
 \end{equation}  
 corresponding to $t$- and $x$-translations and (hyperbolic) rotation,  
 given with respect to the standard pseudo-Cartesian coordinates $(t,x)$.   
Note the generators (\ref{TXH}) of the Lie algebra $i(\mathbb{R}_1^2)$  
 enjoy the following commutator relations:    \begin{equation}  \label{COMM}   
 [{\bf T},{\bf X}] = 0, \quad   [{\bf T}, {\bf H}] = {\bf X}, \quad  
 [{\bf X}, {\bf H}] = {\bf T}.    \end{equation}    
Thus the general form of an element of ${\cal K}^2(\mathbb{R}_1^2)$   
is  given by  \begin{equation} \label{GF}  \begin{array}{rcl}   {\bf K} & = &  
 a_0 {\bf T}\odot {\bf T} + a_1 {\bf T}\odot {\bf X}   + a_2 {\bf X}\odot {\bf X} +  
\\ [0.3cm]    & &   a_3 {\bf T}\odot {\bf H} + a_4 {\bf X} \odot {\bf H}  
+  a_5 {\bf H}\odot {\bf H},    \end{array}   
\end{equation}   
where $\odot$ stands for the symmetric tensor product and   
$a_0, \ldots, a_5 \in \mathbb{R}$ are  arbitrary constants.    
The formula (\ref{GF}) can be used in the problem of  classification   
of the elements of ${\cal K}^2(\mathbb{R}_1^2)$ and thus,  the orthogonal 
 coordinate webs that they generate. For more details,  see Kalnins \cite{Kal}. }  
\end{exa} 
Another aproach that can be used  in the study of Killing tensors of  
valence two is based on  algebraic properties of the matrices that define this type  
of Killing tensors.  Thus, in this case the problem of classification  can be   
 solved by making use of the eigenvalues and eigenvectors of  the Killing tensors; 
 for a complete description of the method  see the review by Benenti \cite{Be} 
 and the related references therein.    These observations provide   
{\em compelling evidence that the study of  Killing tensors lies within 
 the framework of  Riemann's approach to geometry. }  Indeed, Killing tensors 
 appear naturally in Riemann's metric geometry,  as well  as various physical models 
 defined in terms of intrinsic geometry  on pseudo-Riemannian spaces.   
 A new approach  to the study of Killing tensors  introduced in  \cite{MST1} by 
 McLenaghan, Smirnov and The is based on the fact that  Killing tensors of a 
 fixed valence defined on  a pseudo-Riemannian manifold $(M, {\bf g})$  
 of constant curvature  constitute a {\em  vector space}.  
 This  easily follows from the $\mathbb{R}$-bilinear properties of the   
Schouten bracket \cite{Sch} that appears in the fundamental formula  (\ref{KE}).  
Accordingly, one can treat a Killing tensor as an element   of its respective vector space.
     \begin{exa}  {\rm  
Consider again the {\em vector space}  ${\cal K}^2(\mathbb{R}_1^2)$.  
 Solving the Killing tensor equation (\ref{KE})  in the pseudo-Cartesian coordinates  
yields the general formula  \cite{MST3, MST6}  
 \begin{equation}   {\bf K} = \left( \begin{array}{cc}       
a_0 + 2 a_3 x + a_5 x^2         &  a_1 + a_3 t + a_4 x + a_5 tx \\[0.3cm]     
  a_1 + a_3 t + a_4 x + a_5 tx    &  a_2 + 2 a_4 t + a_5 t^2       
 \end{array} \right),  \label{MKT}    \end{equation}   
 of the elements of ${\cal K}^2(\mathbb{R}_1^2)$.  
The arbitrary constants of  integration $a_0,\ldots, a_5$ 
 are the same as in (\ref{GF}), they represent  the dimension   
 of the space ${\cal K}^2(\mathbb{R}_1^2)$.   
The formula (\ref{MKT}) is the ITKT analgoue of the general formula (\ref{2Q}),  
 representing  the elements of the vector space  
 ${\cal Q}^2(\mathbb{R}^2)$ of quadratic forms in CIT.    
} \end{exa}     
We note that in the case of vector spaces of Killing tensors  defined in 
 $\mathbb{R}_1^2$, the generators (\ref{TXH})  are not connected via any non-trivial  
relations.  This is also true for any other two-dimensional pseudo-Riemannian  manifold 
 of constant curvature. In this view,  for a fixed $n\ge 1$ the dimension 
  of the corresponding vector space  ${\cal K}^n(\mathbb{R}_1^2)$  
 can be computed, for example,  by employing the well-known formula  for the 
 dimension of the space $\mbox{Sym}^r(M)$ of symmetric $(r,0)$-tensors  defined  
over an $m$-dimensional manifold:    
\begin{equation} \label{Sym}   
\dim\, \mbox{Sym}^r(M) = {m+r-1 \choose r}.   \end{equation}   
Indeed, in our case  $m = \dim\, i(\mathbb{R}_1^2) = 3$ and $r = n$.   
Therefore we have from (\ref{Sym})  
\begin{equation}  \dim\, {\cal K}^n(\mathbb{R}_1^2) 
 = \frac{1}{2}(n+1)(n+2). \label{dimK}   \end{equation}  
 For spaces of higher dimensions the formula (\ref{dimK}) is no longer valid   
due to the existence of additional  non-trivial relations  among the generators of  
the Lie algebra of Killing vectors  (i.e., the ``syzygy modules problem'' \cite{De}).  
 In the early 1980's the problem of extending the formula (\ref{dimK})  
 to spaces of higher  dimensions was solved independently and   almost  
simultaneously by Delong \cite{De}, Takeuchi \cite{Ta} and Thompson \cite{Th}.  
 According to the {\em Delong-Takeuchi-Thompson (DTT) formula},   
for a fixed $n\ge 1$ the dimension $d$ of the vector space  ${\cal K}^n(M)$ 
 of Killing $(n,0)$-tensors defined  on an $m$-dimensional pseudo-Riemannian 
 manifold $(M, {\bf g})$ is given by   
\begin{equation}   
d = \dim\,{\mathcal K}^n(M)  = \frac{1}{m}{m+n \choose n+1}{m+n-1 \choose n},   
\quad n \ge 1. \label{DTT}    \end{equation}  
 Note the formula (\ref{dimK}) is in agreement with (\ref{DTT}).    
 Having  the  vector spaces of Killing tensors enables one to study them   
under the action of a {\em transformation group}.    
The most natural choice of such a group is, without any doubt,  
the corresponding  Lie group of isometries $I(M)$ of the underlying 
 pseudo-Riemannian manifold  $(M, {\bf g})$.  Indeed,  it is easy to see that 
for a given vector space  ${\cal K}^n(M)$, $n\ge 1$  defined on $(M, {\bf g})$
 the corresponding isometry group   $I(M)$ acts as an automorphism:  
 $I(M):\, {\cal K}^n(M) \rightarrow {\cal K}^n(M)$.  This  key observation  
made by McLenaghan {\em et al} \cite{MST1}  led to the emergence of ITKT.  
More specifically,  the isometry group $I(M)$  acting on $M$  induces the corresponding 
transformation laws on the parameters   $a_0, \ldots, a_{d-1}$  
of the vector space ${\cal K}^n(M)$:   
\begin{equation}  \begin{array}{rcl}   
\tilde{a}_0 & = & \tilde{a}_0(a_0,\ldots, a_{d-1}, g_1, \ldots, g_r),  
\\  \tilde{a}_1 & =  & \tilde{a}_1(a_0,\ldots, a_{d-1}, g_1, \ldots, g_r),  
\\  {} & \vdots & {}  \\  \tilde{a}_{d-1} & = &  \tilde{a}_{d-1}(a_0,\ldots,
 a_{d-1}, g_1, \ldots, g_r).    \end{array} \label{TL}   
\end{equation}   
where $g_1, \ldots, g_r$ are local coordinates on $I(M)$    
that parametrize the group, $r =  \dim\, I(M) = \frac{1}{2}m(m+1)$ and  $d$ 
 is given by (\ref{DTT}).  The formulas (\ref{TL}) can be obtained   in each case by  
employing the standard transformation rules for tensors.    
\begin{exa} {\rm 
Consider again the vector space ${\cal K}^2(\mathbb{R}^2_1)$.    
The corresponding isometry group $I(\mathbb{R}_1^2)$  acts in the Minkowski plane 
 $\mathbb{R}_1^2$ parametrized  by the standard pseudo-Cartesian coordinates  
$(t,x)$ as follows.    \begin{equation}       \left( \begin{array}{c}  \tilde{t}
\\ \tilde{x} \end{array} \right)        = \left( \begin{array}{cc}        
  \cosh\phi & \sinh\phi\\          \sinh\phi & \cosh\phi          
\end{array} \right)           \left( \begin{array}{c} t\\ x\end{array} \right)       
    +  \left( \begin{array}{c} a\\b \end{array}\right),  \label{RM}  
  \end{equation}  
 where $\phi, a, b \in \mathbb{R}$ are local coordinates  that parametrize  
the group $I(\mathbb{R}_1^2)$.  We use the formula (\ref{RM}) and  
 the transformation laws for  $(2,0)$ tensors    
\begin{equation}   
\tilde{K}^{ij}(\tilde{y^1},\tilde{y^2},\tilde{a}_0,\ldots, \tilde{a}_5)   
 = K^{k\ell}(y^1,y^2,a_0,\ldots, a_5)\frac{\partial    
\tilde{y}^i}{\partial y^k}\frac{\partial  \tilde{y}^j}{\partial y^{\ell}},   
\quad i,j, k,\ell = 1, 2, \label{KTR}   \end{equation}   
 where the tensor components  $K^{ij}$ are given by (\ref{MKT}),  
$y^1 = t, y^2 = x$ , to obtain the transformation formulas for the parameters  
 $a_0, \ldots, a_5$ that appear in (\ref{MKT}) and define the parameter space  
 of ${\cal K}^2(\mathbb{R}_1^2)$ \cite{Kal, MST3}:    
\begin{equation}  \label{transf}   \begin{array}{rcl}   
 \tilde{a}_0 & = &  a_0\cosh^2\phi + 2a_1\cosh\phi\sinh\phi + a_2\sinh^2\phi + a_5 b^2 
 \\[0.3cm]    & & -2 (a_3\cosh \phi  + a_4\sinh \phi)b, 
\\[0.3cm]      \tilde{a}_1 &=&  a_1(\cosh^2\phi  + \sinh^2\phi)  
 +(a_0 + a_2)\cosh\phi\sinh\phi \\[0.3cm]    & &  
- (a a_3 + b a_4)\cosh\phi - (a a_4 + b a_3)\sinh\phi + a_5 ab,
\\[0.3cm]      \tilde{a}_2 & = &  a_0\sinh^2\phi + 2a_1\cosh\phi\sinh\phi 
 + a_2\cosh^2\phi + a_5 a^2 \\[0.3cm]    & & -2 (a_4\cosh \phi  + a_3\sinh \phi)a,  
\\[0.3cm]       \tilde{a}_3 & = &  a_3\cosh\phi + a_4\sinh\phi - a_5 b,   
\\[0.3cm]   \tilde{a}_4 & = & a_3\sinh \phi + a_4 \cosh\phi -a_5 a,   
\\[0.3cm]   \tilde{a}_5 &=& a_5.    \end{array}   
\end{equation}   
} \end{exa}    We note that the corresponding transformation formulas for 
the parameters obtained in \cite{MST3}  were derived for  
{\em covariant} Killing tensors. Accordingly, they differ somewhat from  (\ref{transf}) 
presented above (compare with (7.6) in \cite{MST3}). Clearly, the transformation
 formulas (\ref{transf}), and   more generally - (\ref{TL}) are analogues of 
the corresponding transformation   formulas in CIT (see, for example, (\ref{TLP})). 
It must be mentioned, however,  that in the  case of ITKT they are  
 computationally more difficult to obtain.    In view of the above observations,
 it is now easy to determine the ITKT analogue   of the CIT-concept of an invariant.   
 \begin{defi} Let $(M,{\bf g})$  be a  pseudo-Riemannian manifold 
 of constant curvature.  For a fixed $n \ge 1$ consider  the corresponding  
 space ${\cal K}^n(M)$ of Killing tensors of valence $n$   defined in $(M,{\bf g})$. 
A smooth function   ${\cal I}:\,  \Sigma \rightarrow \mathbb{R}$  defined in
 the space of functions  over the parameter space  $\Sigma$ is said to be an 
$I(M)$-{\em invariant  of the vector space  ${\cal K}^n(M)$} iff it satisfies the 
condition   
\begin{equation}   
{\cal I} = F(a_0, \ldots, a_{d-1})   
= F(\tilde{a}_o, \ldots, \tilde{a}_{d-1}) \label{Inv}   
\end{equation}   
 under the transformation laws (\ref{TL}) induced by the isometry group $I(M)$.    
\end{defi}   
 We note that in a similar way the ITKT-analogues of the CIT-concepts   of a 
{\em covariant} and {\em joint invariant}  have been introduced in \cite{SY}.     
In  complete analogy with CIT, we can in principle determine  the space  of
 $I(M)$-invariants for a specific vector space  ${\cal K}^n(M)$, $n\ge 1$ by  
employing the (Sophus Lie) method of  infinitesimal generators. To do so, one has  
to determine the (infinitesimal)  action of $I(M)$ in the corresponding parameter  space 
 $\Sigma$ of ${\cal K}^n(M)$  defined by the parameters   $a_0, \ldots, a_{d-1}$.
 McLenaghan {\em et al} \cite{MST1}  
devised an original  procedure  that can  be used to derive the generators 
 of the Lie algebra in  $\Sigma$ isomorphic to the Lie algebra $i(M)$ of  $I(M)$ 
and thus,  to compute  the invariants. We briefly review  the {\em MST-procedure} 
\cite{MST1} here.   Let $\{{\bf X}_1, \dots, {\bf X}_r\}$  be the  infinitesimal 
generators  (Killing vector fields)  of the Lie group $I(M)$  acting on $M$. Note  
 $\mbox{Span}\, \{ {\bf X}_1, \ldots, {\bf X}_r\} = {\cal K}^1(M) = i(M)$,  
 where $i(M)$ is the Lie algebra of the Lie group $I(M)$. For a fixed $n\ge 1$,  
 consider the corresponding vector space ${\cal K}^n(M)$.  To determine  the action  
of $I(M)$ in the parameter space $\Sigma$ defined  by $a_0, \ldots, a_{d-1}$,  
we find first the  infinitesimal generators of  $I(M)$ in $\Sigma$. 
 Consider $\Diff\, \Sigma$, it defines the corresponding space   $\Diff\,
 {\cal K}^n(M)$, whose elements are determined by the elements of   $\Diff\, \Sigma$ 
in an obvious way. Let ${\bf K}^0 \in \Diff\, {\cal K}^n(M)$.   Note ${\bf K}^0$ is 
determined by $d$ parameters $$a_i^0= a_i^0(a_0,\ldots, a_{d-1}),$$  where 
 $i=0,\ldots, d-1$, which are functions of  $a_0, \ldots, a_{d-1}$ - the  parameters 
  of $\Sigma$. Define now a map  $\pi:\, \Diff \, 
{\cal K}^n(M) \rightarrow {\cal X}(\Sigma),$    given by    
\begin{equation}   
{\bf K}^0 \rightarrow \sum_{i=0}^{d-1}a_i^0(a_0,\ldots, a_{d-1}) {\partial}_{a_i}. 
 \label{pi}   \end{equation}   
To specify the action of $I(M)$ in $\Sigma$, we have to find the  counterparts 
 of the generators  ${\bf X}_1, \ldots, {\bf X}_r$ in ${\cal X}(\Sigma)$.  
 Consider the composition $\pi\circ {\cal L}$, where $\pi$ is defined  
 by (\ref{pi}) and ${\cal L}$ is the Lie derivative operator.  
 Let  ${\bf K}$ be the general Killing tensor of ${\cal K}^n(M)$,  
 in other words $\bf K$ is the general solution to  the Killing tensor equation 
 (\ref{KE}).   Next, define   \begin{equation} \label{comv1}   
{\bf V}_i = \pi{\cal L}_{{\bf X}_i}\,{\bf K}, \quad i = 1, \ldots r.    
\end{equation}   
The composition map   
\begin{equation} \label{compm} 
 \pi \circ {\cal L}:\, i(M) \rightarrow {\cal X}(\Sigma)  
\end{equation}   maps the generators ${\bf X}_1, \ldots, {\bf X}_r$ to 
${\cal X}(\Sigma)$.  Finally, we check that the vector fields ${\bf V}_i,$ $i=1, 
\ldots, r$   satisfy the same commutator relations as the original  
 ${\bf X}_i,$ $i=1, \ldots, r$. This step is actually redundant,   
since it has been proven in general by showing that Killing tensors   
can be expressed as irreducible representations of $GL(n, \mathbb{R})$    
that the vector fields (\ref{comv1}) satisfy the same commutator relations  
 as the original generators of $i(M)$ \cite{MMS}. The main result of \cite{MMS}  
 is the proof of the corresponding conjecture formulated in \cite{MST3}.  
 Therefore we can use the vector fields (\ref{comv1}) to solve the problem   
 of the determination of the $I(M)$-invariants of the vectors space   
${\cal K}^n(M)$ employing the (Sophus Lie)   {\em method of infinitesimal generators}
 by solving the   corresponding system of (linear) PDEs generated by (\ref{comv1}): 
  \begin{equation} \label{system}   
{\bf V}_i (F) = 0, \quad i = 1, \ldots, r.   \end{equation}   
 The general solution to the system (\ref{system}) describes the space  
 of all $I(M)$-invariants of the vector space ${\cal K}^n(M)$, $n\ge 1$.   
 The MST-procedure \cite{MST1} based on {\em  Lie derivative deformations of Killing
 tensors}   is an analogue of the standard exponentiation used in CIT 
to determine  the action  of a group in the parameter  space determined by the vector 
 space of  homogeneous polynomials in question. The technique of the  Lie derivative 
 deformations used here is a very powerful tool.  It was used before, for example, 
 in Smirnov \cite{Sm}  to generate compatible Poisson bi-vectors in  the theory of
  bi-Hamiltonian systems. The idea introduced in \cite{Sm}  was used 
 in  \cite{Ser} and applied to a different class of  integrable systems.   
We also note that  the generators (\ref{system})  can  be alternatively determined 
 from the parameter  transformation laws (\ref{TL}) when the latter are available.  
 However, it is  increasingly  difficult and often impossible to  
 determine (\ref{comv1}) using (\ref{TL}) for  vector spaces of  Killing tensors 
 of higher valence or defined in pseudo-Riemannian  manifolds of higher dimensions. 
 In what follows, we  employ the  MST-procedure to prove the ITKT-analogue of  
Cayley's lemma that  is Problem \ref{Prob}  formulated in the previous section. 
 To illustrate the effectiveness  of the MST-procedure,   let us  consider the 
following example.    \begin{exa} {\rm  Consider again the vector space 
${\cal K}^2(\mathbb{R}_1^2)$.   The action of the isometry group $I(\mathbb{R}_1^2)$
 on the corresponding   parameter space $\Sigma$ defined by $a_0, \ldots, a_5$ 
(see (\ref{MKT}))   is determined by the formulas (\ref{transf}).  In order to determine 
the  infinitesimal action of $I(\mathbb{R}_1^2)$  in $\Sigma$, we employ  
the MST-procedure.   Thus, using the general formula  (\ref{MKT}) in conjunction
 with (\ref{comv1}),   we derive the corresponding generators ${\bf V}_i,$ $i=1,2,3$:    
\begin{equation} \label{vecM}   \begin{array}{rcl}   
{\bf V}_1 &=&  a_3\partial _{a_1}+2a_4\partial_{a_2}+a_5 \partial_{a_4},  
\\[0.3cm]   {\bf V}_2 &=&  a_4\partial _{a_1}+2a_3\partial_{a_0}+a_5 \partial_{a_3},
 \\[0.3cm]   {\bf V}_3 &=&  -2a_1\partial _{a_0}-a_4\partial _{a_3} 
   -(a_0+a_2)\partial _{a_1}-2a_1 \partial _{a_2}-a_3\partial _{a_4}.   
\end{array}   \end{equation}    
We immediately note that the vector fields (\ref{vecM})  satisfy the  following 
commutator relations:   
$$[{\bf V}_1, {\bf V}_2] = 0, \quad [{\bf V}_1, {\bf V}_3] = - {\bf V}_2,  
 \quad [{\bf V}_2, {\bf V}_3] = - {\bf V}_1.$$  
 Choosing the basis $\{ - {\bf V}_1, - {\bf V}_2, - {\bf V}_3\}$  
 reveals that the Lie algebra generated by (\ref{vecM}) is isomorphic to   
the Lie algebra $i(\mathbb{R}_1^2)= {\cal K}^1(\mathbb{R}_1^2)$ generated  
 by (\ref{TXH}). Indeed, the vector fields  (\ref{TXH}) and 
 $\{ - {\bf V}_1, - {\bf V}_2, - {\bf V}_3\}$,   where ${\bf V}_i,$ 
$i = 1,2,3$ are given by (\ref{vecM}) satisfy the same   commutator relations
 (see (\ref{COMM})). We conclude therefore that    the vector fields (\ref{vecM})
 represent the infinitesimal action  of  $I(\mathbb{R}_1^2)$ in $\Sigma$. 
Our next observation is that in view  of  Proposition \ref{Prop1} the orbits of the
 $I(\mathbb{R}_1^2)$-action   have dimension three wherever the vector fields 
(\ref{vecM})  are linearly independent.  Therefore in that subspace of $\Sigma$, 
 by Theorem \ref{FT}, we expect  to derive   6-3 = 3 fundamental 
 $I(\mathbb{R}_1^2)$-invariants.  The infinitesimal generators  of the  
$I(\mathbb{R}_1^2)$-action in the 5-dimensional vector subspace of   
{\em non-trivial Killing two tensors}  of ${\cal K}^2(\mathbb{R}_1^2)$  
 were determined in McLenaghan {\em et al} \cite{MST4, MST6}.     
} \end{exa}  
 Employing the method of characteristics to solve the system of PDEs (\ref{system}) 
 defined by the vector fields  (\ref{vecM}), we arrive  at the following theorem.    
\begin{theo}  
 Any algebraic $I(\mathbb{R}^2_1)$-invariant ${\cal I}$ of the subspace  
 of the parameter space $\Sigma$ of ${\cal K}^2(\mathbb{R}_1^2)$  defined   
by the condition that the vector fields (\ref{vecM})  are linearly independent  
 can be (locally) uniquely expressed  as an analytic function   
$${\cal I}= F(\Delta_1, \Delta_2, \Delta_3)$$  where the fundamental invariants 
 $\Delta_i,$ $i=1,2,3$  are given by   
\begin{equation}  \begin{array}{rcl}   
\Delta_1 & = & a_5, 
\\ [0.3cm]   \Delta_2 & = & (a_0-a_2)a_5-a_3^2+a_4^2,   
\\ [0.3cm]   \Delta_3 & = & (a_3^2+a_4^2-a_5(a_0+a_2))^2 -4(a_5a_1-a_3a_4)^2.   
\end{array}  \end{equation}    
\end{theo}    
The fact that $\Delta_1= a_5$ is a fundamental  $I(\mathbb{R}_1^2)$-invariant  
of the vector space  ${\cal K}^2(\mathbb{R}_1^2)$ can be trivially deduced from 
  the transformation formulas (\ref{transf}). The fundamental   
 $I(\mathbb{R}_1^2)$-invariant $\Delta_3$ presented above was first  derived 
 in McLenaghan {\em et al} \cite{MST4, MST6} and used to  generate {\em discrete}  
$I(\mathbb{R}_1^2)$-invariants,  which were in turn employed to classify orthogonal  
coordinate webs  in the Minkowski plane $\mathbb{R}_1^2$.  
 The same problem was  solved in \cite{SY} by employing the 
 $I(\mathbb{R}_1^2)$-invariants  and covariants of the vector space 
 ${\cal K}^2(\mathbb{R}_1^2)$.    The observations and results summarized above 
 put in evidence  that  {\em ITKT is a part of F. Klein's approach to geometry}.   
 This is especially evident when one considers the vector spaces of Killing tensors  
 of valence two.  Thus, for example, in Horwood {\em et al} \cite{HMS} 
 the orthogonal  coordinate webs of the Euclidean space $\mathbb{R}^3$  
were completely classified  in terms of the $I(\mathbb{R}^3)$-invariants. 
 This is something to be expected  since the theory of orthogonal  coordinate  
webs of $\mathbb{R}^3$ is a part of  the Euclidean geometry which,  
according to  Felix Klein's ``Erlangen Program''  \cite{FK72, FK93}, 
is an invariant theory of the  corresponding isometry group  $I(\mathbb{R}^3)$.  
 In Section 3 we use the  results presented above,  in particular, the MST-procedure, 
 to solve Problem  \ref{Prob}.      
 \section{The ITKT analogue of Cayley's lemma}     
 In this section we prove the ITKT analogue of the Cayley Lemma \cite{AC} 
 presented in Section 1. The vector space ${\cal K}^n(\mathbb{R}_1^2)$ appears  
to be a natural counterpart of the vector space ${\cal Q}^n(\mathbb{R}^2)$ in CIT. 
 The problem can be solved by employing  the MST-procedure described in the  previous
 section.  To proceed, we need to derive first a general formula for the  elements of
 ${\cal K}^n(\mathbb{R}_1^2)$ (i.e., an analogue of (\ref{gQ})).  
Note, that by (\ref{dimK}) the dimension of the vector space in question is 
 $(n+1)(n+2)/2$. Thus, each contravariant tensor  
${\bf K} \in {\cal K}^n(\mathbb{R}_1^2)$ is determined by $(n+1)(n+2)/2$
  parameters that appear in the $n+1$  components of the form   
\begin{equation} \label{comp} K^{i_1\ldots i_p j_1\ldots j_{n-p}},   
\end{equation}  where $i_1 =\cdots = i_p = 1$,  $j_1 = \cdots = j_{n-p} = 2$ 
and $p = 0,1, \ldots, n$.  To derive the formulas for the components (\ref{comp}), 
we solve the Killing  tensor equation (\ref{KE}) in the coordinates $(t, x)$,
 which in this case reduces  to the following system of PDEs:  
 \begin{equation}  \left\{\begin{array}{l}   
\partial _t{K}^{i_1\cdots i_n} = 0,   
\\ [0.3cm]  (n-p+1)\partial_x {K}^{i_1\cdots i_{p}j_1\cdots j_{n-p}}  
=   p\partial_t {K}^{i_1\cdots i_{p-1} j_1\cdots j_{n-p+1}},  
\\ [0.3cm] \partial _x{K}^{j_1\cdots j_n} = 0,  
\end{array}\right.  \label{PDES} \end{equation}   
where $p = 0, 1, \ldots, n$, $\partial _t=\frac{\partial}{\partial t},$  
 $\partial _x=\frac{\partial}{\partial _x}.$ As a consequence of (\ref{PDES}), 
we  readily obtain the necessary differential conditions: 
\begin{equation} \begin{array}{l}
 (\partial_x)^{p+1}K^{i_1\ldots i_p j_1\ldots j_{n-p}} = 0, 
\\[0.3cm] (\partial_t)^{n-p+1}K^{i_1\ldots i_p j_1\ldots j_{n-p}} =0. 
\end{array} \label{NC} \end{equation}  
Solving (\ref{NC}), we arrive  at the following result. Each component of 
 (\ref{comp}) is a mixed polynomial of degree $p$ in $x$ and degree $q$ 
in $t$: :
\begin{equation}   K^{i_1\cdots i_p j_1\cdots j_q}   =
 \left\{ \begin{array}{l @{\quad \mbox{if}\quad}l}   \D \sum_{i=0}^q {p\choose i} t^i  
 \sum_{j=0}^p{p\choose j}a_{pij}x^j,   & \displaystyle n\ge p\ge \Big[\frac{n+1}{2}\Big], 
 \\[0.6cm]  \D \sum_{i=0}^p {q\choose i} x^i   \sum_{j=0}^q {q\choose j}b _{qij}t^j,  
 & \displaystyle 0\le p \le  \Big[\frac{n+1}{2}\Big],   \end{array}\right.  
\label{formula}  
\end{equation}    
 where $q = n - p$ and  the parameters $a_{pij}, b_{qij}$ are to be determined 
(at this stage they are inserted for mere convenience).  
 We immediately recognize that the formula (\ref{formula}) is the ITKT analogue of 
 the general formula (\ref{gQ}) exhibited in Section 1.  
The parameters $a_{pij}, b_{qij}$ can be determined by following the general 
 procedure of solving the system of PDEs (\ref{PDES}). For convenience we consider 
separately two cases: $n = 2k + 1$ and $n = 2k$. The parameters of each of the $n+1$ 
components can be  organized into groups in such a way that the parameters of 
  one group are completely  determined by the parameters of the other 
(see the illustrative examples below).   After relabelling the parameters,
 we arrive at the following two schemes  (corresponding to $n = 2k$ and $n = 2k+1$ 
respectively), which specify the arrangements of the parameters of the first groups
 of the components.  Once they are specified, the parameters of the other groups can
 be determined accordingly.   
\bigskip  \textbf{Case 1:} $\quad n=2k$  \begin{equation} \label{C2}   
\begin{array}{l@{\qquad }l@{\quad }l@{\quad }l@{\quad }l@{\quad }l@{\quad }l}   
\mbox{Step }1: & a^1_0 & a^1_1 &\cdots & a^1_{n-2} & a^1_{n-1}& a^1_n,  
\\[0.3cm] {} & b^1_0 & b^1_1 & \cdots & b^1_{n-2} & b^1_{n-1} & a^1_n   
\\[0.6cm]  \mbox{Step }2: & a^2_0 & a^2_1 & \cdots & a^2_{n-3} & a^2_{n-2} & b^1_{n-1}   
\\[0.3cm] {} & b^2_0 & b^2_1 & \cdots  & b^2_{n-3}& a^2_{n-2}& a^1_{n-1}   
\\[0.6cm]  {} & {} &{} & \vdots & {} & {} & {}   
\\[0.6cm]\D  \mbox{Step } \frac{n}{2}:  & a^{\frac{n}{2}}_0 & a^{\frac{n}{2}}_1  
 & a^{\frac{n}{2}}_2 & b^{\frac{n-2}{2}}_1  & \cdots & b^1_{\frac{n+2}{2}}   
\\[0.3cm] {}& b^{\frac{n}{2}}_0 & b^{\frac{n}{2}}_1 & a^{\frac{n}{2}}_2  
& a^{\frac{n-2}{2}}_1 & \cdots & a^1_{\frac{n+2}{2}}   
\\[0.6cm]\D  \mbox{Step} \frac{n+2}{2}: & a^{\frac{n+2}{2}}_0 
 & b^{\frac{n}{2}}_1 & b^{\frac{n-2}{2}}_2 & b^{\frac{n-4}{2}}_3 
 & \cdots & b^1_{\frac{n}{2}}   
\end{array}  \end{equation}     
\textbf{Case 2:} $\quad n=2k+1$  
\begin{equation}  \label{C1}  
\begin{array}{l@{\qquad }l@{\quad }l@{\quad }l@{\quad }l@{\quad }l@{\quad }l}    
\mbox{Step 1}: &  a^1_0 & a^1_1 & \cdots  & a^1_{n-2} & a^1_{n-1} & a^1_n,   
\\[0.3cm]{} & b^1_0 & b^1_1 & \cdots  & b^1_{n-2} & b^1_{n-1} & a^1_n   
\\[0.6cm]  \mbox{Step }2: &  a^2_0 & a^2_1 & \cdots & a^2_{n-3} & a^2_{n-2}& b^1_{n-1}   
\\[0.3cm] {} & b^2_0 & b^2_1 & \cdots & b^2_{n-3} & a^2_{n-2} & a^1_{n-1}   
\\[0.6cm]  {} & {} & {}& \vdots & {} & {} & {}  
 \\[0.6cm] \D \mbox{Step} \frac{n-1}{2}:   &  a^{\frac{n-1}{2}}_0
 & a^{\frac{n-1}{2}}_1  & a^{\frac{n-1}{2}}_2  & a^{\frac{n-1}{2}}_3 
 & \cdots & b^1_{\frac{n+3}{2}}   
\\[0.3cm] {} & b^{\frac{n-1}{2}}_0 & b^{\frac{n-1}{2}}_1  
& b^{\frac{n-1}{2}}_2 &  a^{\frac{n-1}{2}}_3  & \cdots &  a^1_{\frac{n+3}{2}}   
\\[0.6cm]\D  \mbox{Step}\frac{n+1}{2}:  & a^{\frac{n+1}{2}}_0 
 & a^{\frac{n+1}{2}}_1  & b^{\frac{n-1}{2}}_2 & b^{\frac{n-3}{2}}_3  
 & \cdots & b^1_{\frac{n+1}{2}}  
 \\[0.3cm]\D {}& b^{\frac{n+1}{2}}_0 & a^{\frac{n+1}{2}}_1  
& a^{\frac{n-1}{2}}_2 & a^{\frac{n-3}{2}}_3 & \cdots & a^1_{\frac{n+1}{2}}   
\end{array}  \end{equation}         
The parameters that appear in the general solution to (\ref{formula}) 
 are now organized in two schemes according the cases of $n$ being even (\ref{C1})  
 and odd (\ref{C2}) respectively. More specifically, 
we first give  $2(n+1)-1$ parameters 
$$a^1_0,\ldots,a^1_{n-1},a^1_n b^1_0,\ldots,b^1_{n-1},a^1_n,$$
 and then  write down the first and the last components of the general element 
${\bf K}\in {\cal K}^n({\mathbb R}^2_1)$ as follows: 
\begin{eqnarray*} 
K^{11\cdots 11}&=& \left[a^1_0+{n\choose 1}a^1_1x 
+{n\choose 2}a^1_2x^2+\ldots +{n\choose n-1}a^1_{n-1}x^{n-1}+a^1_nx^n\right],
\\[0.4cm] K^{22\cdots 22}&=& \left[b^1_0+{n\choose 1}b^1_1t+{n\choose 2}b^1_2t^2
+\ldots +{n\choose n-1}b^1_{n-1}t^{n-1}+a^1_nt^n\right]. 
\end{eqnarray*} 
Next step: For  $2(n-1)-1$ new parameters
 $$a^2_0,\ldots,a^2_{n-3}, a^2_{n-2}, b^2_0,\ldots,b^2_{n-3}, a^2_{n-2}$$ 
 we  then write down the second  and penultimate components of 
 ${\bf K}$ as follows (see (\ref{formula})),
 each of which is the sum  of two   polynomials, the first having
 been determined by the newly specified parameters and the other  - by the
 parameters determined previously.  
\begin{eqnarray*} K^{11\cdots 12}&=& \left[a^2_0+{n-1\choose 1}a^2_1x +\ldots 
+{n-1\choose n-2}a^2_{n-2}x^{n-2}+b^1_{n-1}x^{n-1}\right]\nonumber 
\\[0.4cm] & & +t\left[a^1_1+{n-1\choose 1}a^1_2x +\ldots 
+{n-1\choose n-2}a^1_{n-1}x^{n-2}+a^1_{n}x^{n-1}\right],
\\[0.4cm] K^{22\cdots 21}&=& \left[b^2_0+{n-1\choose 1}b^2_1t +\ldots 
+{n-1\choose n-2}a^2_{n-2}t^{n-2}+a^1_{n-1}t^{n-1}\right]\nonumber 
\\[0.4cm] && +x\left[b^1_1+{n-1\choose 1}b^1_2t +\ldots 
+{n-1\choose n-2}b^1_{n-1}t^{n-2}+a^1_{n}t^{n-1}\right]. 
\end{eqnarray*}  To clarify the process more, let us consider  the next step 
(if any)for the given $2(n-3)-1$ parameters 
\[a^3_0,a^3_1,\ldots,a^3_{n-5},a^3_{n-4}, b^3_0,b^3_1,\ldots,b^3_{n-5},a^3_{n-4} \] 
we write down the next two comonents as follows: 
\begin{eqnarray*} \lefteqn{ K^{11\cdots 122}=}
\\ && \left[a^3_0+{n-2\choose 1}a^3_1x +\ldots +{n-2\choose n-4}a^3_{n-4}x^{n-4} 
+{n-2\choose n-3}b^2_{n-3}x^{n-3} +                b^1_{n-2}x^{n-2}\right]
\\[0.4cm]&& +2t\left[a^2_1+{n-2\choose 1}a^2_2x 
+\ldots +{n-2\choose n-4}a^2_{n-3}x^{n-4} +{n-2\choose n-3}a^2_{n-2}x^{n-3} +                b^1_{n-1}x^{n-2}\right]\\[0.4cm]&& +t^2\left[a^1_2+{n-2\choose 1}a^1_3x +\ldots +{n-2\choose n-4}a^1_{n-2}x^{n-4} +{n-2\choose n-3}a^1_{n-1}x^{n-3} +                a^1_{n}x^{n-2}\right], \end{eqnarray*} \begin{eqnarray*} \lefteqn{ K^{22\cdots 211}=}\\ && \left[b^3_0+{n-2\choose 1}b^3_1t +\ldots +{n-2\choose n-4}a^3_{n-4}t^{n-4} +{n-2\choose n-3}a^2_{n-3}t^{n-3} +                a^1_{n-2}t^{n-2}\right]\\[0.4cm]&& +2x\left[b^2_1+{n-2\choose 1}b^2_2t +\ldots +{n-2\choose n-4}b^2_{n-3}t^{n-4} +{n-2\choose n-3}a^2_{n-2}t^{n-3} +                a^1_{n-1}t^{n-2}\right]\\[0.4cm]&& +x^2\left[b^1_2+{n-2\choose 1}b^1_3t +\ldots +{n-2\choose n-4}b^1_{n-2}t^{n-4} +{n-2\choose n-3}b^1_{n-1}t^{n-3} +                a^1_{n}t^{n-2}\right]. \end{eqnarray*}          We repeat this process in both  directions (i.e., going ``downwards'' and ``upwards'') until  it is terminated in the middle of (\ref{formula}). In this view, counting the steps in both cases ($n$ is even and $n$ is odd), it is easy to see that the dimension of the space $$d = {\cal K}^n({\mathbb R}^2_1)  = \frac{1}{2}(n+1)(n+2), \quad n \ge 1$$  gets decomposed as follows.  \begin{equation} d=\left\{ \begin{array}{l@{\quad }l} [2(n+1)-1]+[2(n-1)-1]+\cdots +[2\times  2-1] & \mbox{if $n$ is odd,}\\[0.3cm]  [2(n+1)-1]+[2(n-1)-1]+\cdots +[2\times  1-1] & \mbox{if $n$ is even.} \end{array} \right. \end{equation} The auxiliary problem of  finding the general form for the elements ${\bf K}\in {\cal K}^n({\mathbb R}^2_1)$ is therefore completely solved.        We immediately notice that the coefficients in the  general solution (\ref{MKT}) can be relabeled following the scheme (\ref{C2})  as follows: $a_0 = a^1_0$, $a_1=a^2_0$,  $a_2 = b_0^1$, $a_3 = a_1^1$, $a_4 = b_1^1$ and  $a_5 = a_2^1$. To illustrate our results, let us consider more challenging examples.       \begin{exa}{\rm   Consider the vector space ${\cal K}^4({\mathbb R}^2_1)$, note     $d= \dim\, {\cal K}^4({\mathbb R}^2_1)= (4+1)(4+2)/2=15$. Following the coefficient  scheme (\ref{C2}), we arrive at the following  formulas for the components  of the elements of  ${\cal K}^4({\mathbb R}^2_1)$.   \begin{equation}  \label{K4} \begin{array}{rcl}   K^{1111} & = &  a^1_0+4a^1_1x+6a^1_2x^2+4a^1_3x^3+a^1_4x^4,   \\ [0.3cm]  K^{1112} & = &  (a^2_0+3a^2_1x+3a^2_2x^2+b^1_3x^3)+t(a^1_1  +3a^1_2x+3a^1_3x^2+a^1_4x^3),   \\ [0.3cm]  {K}^{1122} & = & (a^3_0+2b^2_1x+b^1_2x^2)+2t(a^2_1+2a^2_2x+b^1_3x^2)  \\ [0.3cm]   & & +t^2(a^1_2+2a^1_3x+a^1_4x^2),   \\ [0.3cm]  K^{1222}  & = & (b^2_0+3b^2_1t+3a^2_2t^2+a^1_3t^3)+x(b^1_1  +3b^1_2t+3b^1_3t^2+a^1_4t^3),   \\ [0.3cm] K^{2222}  & = &  b^1_0+4b^1_1t+6b^1_2t^2+4b^1_3t^3+a^1_4t^4.   \end{array}  \end{equation}   }\end{exa}       \begin{exa}{\rm  Consider the vector space ${\cal K}^5(\mathbb{R}^2_1)$. In this case  $d= \dim \, {\cal K}^5(\mathbb{R}^2_1)=  (5+1)(5+2)/2=21$ and the components are given by    \begin{equation}  \label{K5} \begin{array}{rcl}   K^{11111} & = &  a^1_0+5a^1_1x+10a^1_2x^2+10a^1_3x^3+5a^1_4x^4+a^1_5x^5,   
\\ [0.3cm]  K^{11112} & = &     (a^2_0+4a^2_1x+6a^2_2x^2+4a^2_3x^3+b^1_4x^4) 
\\[.3cm] && +t(a^1_1 +4a^1_2x+6a^1_3x^2+4a^1_4x^3+a^1_5x^4),   
\\ [0.3cm]  {K}^{11122} & = &  
(a^3_0+3a^3_1x+3b^2_2x^2+b^1_3x^3)+2t(a^2_1+3a^2_2x+3a^2_3x^2+b^1_4x^3)  
\\ [0.3cm] &&+t^2(a^1_2+3a^1_3x+3a^1_4x^2+a^1_5x^3),  
\\[0.3cm]  {K}^{11222} & = &  
(b^3_0+3a^3_1t+3a^2_2t^2+a^1_3t^3)+2x(b^2_1+3b^2_2t+3a^2_3t^2+a^1_4x^3)  
\\ [0.3cm] &&+x^2(b^1_2+3b^1_3t+3b^1_4t^2+a^1_5t^3),  
\\[0.3cm]   K^{12222} & = &      
(b^2_0+4b^2_1t+6b^2_2t^2+4a^2_3t^3+a^1_4t^4) 
\\[.3cm] &&  +x(b^1_1 +4b^1_2t+6b^1_3t^2+4b^1_4t^3+a^1_5t^4),   
\\ [0.3cm]  K^{22222} & = & 
 b^1_0+5b^1_1t+10b^1_2t^2+10b^1_3t^3+5b^1_4t^4+a^1_5t^5.   
\end{array}  \end{equation}      
} \end{exa}     
In principle, based on the formulas (\ref{formula}), (\ref{C2}) and (\ref{C1}),
we can now write down  explicitly the general form of the elements of  
${\cal K}^n(\mathbb{R}_1^2)$ for an arbitrary $n$, without any difficulty,  
just following the parameter scheme given above.    To solve  Problem \ref{Prob}, 
 we employ the  MST-procedure \cite{MST1} outlined in the previous section.
 Using the formulas  (\ref{comv1}), (\ref{formula}), (\ref{C2}) and (\ref{C1}), 
we arrive at the  general formulas for the vector fields representing the 
infinitesimal action of  the isometry group $I(\mathbb{R}_1^2)$ 
on the parameter space. As above, we have  two cases corresponding to 
(\ref{C2}) and (\ref{C1})  respectively.     
\medskip  
{\bf Case 1}$\quad n=2k$  
\begin{equation} \label{V1e}  \begin{array}{rcl}  
{\bf V}_{1} & = &  
a^1_1\partial _{ a^2_0} +a^1_2\partial_{ a^2_1} +\cdots +a^1_{n-1}\partial_{  a^2_{n-2}}  
 \\[0.3cm]  & &  \D +2a^2_1\partial _{ a^3_0} +2a^2_2\partial _{ a^3_1} 
 +\cdots +2a^2_{n-3}\partial _{ a^3_{n-4}}   
\\[0.3cm] & & \dotfill    
\\[0.3cm]  & & \D +\frac{n}{2}a^{\frac{n}{2}}_1\partial _{ a^{\frac{n+2}{2}}_0}  
\\[0.3cm]  & & \D +\frac{n+2}{2}b^{\frac{n}{2}}_1\partial _{ b^{\frac{n}{2}}_0} 
 +\frac{n}{2}a^{\frac{n}{2}}_2\partial _{ b^{\frac{n}{2}}_1}   
\\[0.3cm]  & & \dotfill    \\[0.3cm]  & &  \D  +(n-1)b^2_1\partial _{ b^2_0}
 +(n-2)b^2_2\partial _{ b^2_1} +\cdots 2b^2_{n-2}\partial _{ b^2_{n-3}}   
\\[0.3cm]  && \D +nb^1_1\partial _{ b^1_0} +(n-1)b^1_2\partial _{ b^1_1} 
+\cdots +a^1_n\partial _{ b^1_{n-1}},    
\end{array}  \end{equation}     
\begin{equation} \label{V2e}  \begin{array}{rcl}  
{\bf V}_{2} & = &   
 b^1_1\partial _{ b^2_0} +b^1_2\partial _{ b^2_1} 
+\cdots +b^1_{n-1}\partial _{ a^2_{n-2}}   
\\[0.3cm]  &&  \D  +2b^2_1\partial _{\D b^3_0} +2b^2_2\partial _{ b^3_1} 
 +\cdots +2b^2_{n-3}\partial _{ a^3_{n-4}}   
\\[0.3cm] & & \dotfill   
 \\[0.3cm] & & \D +\frac{n}{2}b^{\frac{n}{2}}_1\partial _{ a^{\frac{n+2}{2}}_0}   
\\[0.3cm] && \D +\frac{n+2}{2}a^{\frac{n}{2}}_1\partial _{ a^{\frac{n}{2}}_0} 
 +\frac{n}{2}a^{\frac{n}{2}}_2\partial _{ a^{\frac{n}{2}}_1}   
\\[0.3cm] && \dotfill   
\\[0.3cm] && \D +(n-1)a^2_1\partial _{ a^2_0} +(n-2)a^2_2\partial _{ a^2_1}
 +\cdots +2a^2_{n-2}\partial _{ a^2_{n-3}}   
\\[0.3cm]  && \D +na^1_1\partial _{ a^1_0} +(n-1)a^1_2\partial _{ a^1_1} 
+\cdots +a^1_n\partial _{ a^1_{n-1}},  
\end{array}  \end{equation}   \begin{equation} \label{V3e}  
\begin{array}{rcl}  
{\bf V}_{3} &  = &   
-na^2_0\partial _{ a^1_0} -(n-1)a^2_1\partial _{ a^1_1}-\cdots 
 -2a^2_{n-2}\partial _{ a^1_{n-2}} -b^1_{n-1}\partial _{ a^1_{n-1}}   
\\[0.3cm] && \D  -[(n-1)a^3_0+a^1_0 ]\partial _{ a^2_0}   
 -[(n-2)a^3_1+a^1_1 ]\partial _{ a^2_1}   
\\ && \D  - \cdots -[2b^2_{n-3}+a^1_{n-3} ]\partial _{ a^2_{n-3}}  
 \\[0.3cm] & & \dotfill   \\[0.3cm]   & & \D  -\frac{n}{2} [a^{\frac{n}{2}}_0
 +b^{\frac{n}{2}}_0 ] \partial _{ a^{\frac{n+2}{2}}_0}-\cdots 
 -[a^1_{n-2}+b^1_{n-2}]\partial _{ a^2_{n-2}}   
\\[0.3cm] & & \dotfill    
\\[0.3cm] & & \D  -[(n-1)b^3_0+b^1_0 ]\partial _{b^2_0}    
 - \cdots -[2a^2_{n-3}+b^1_{n-3} ]\partial _{ b^2_{n-3}}   
\\[0.3cm] && \D  -nb^2_0 \partial _{ b^1_0} -(n-1)b^2_1 \partial _{ b^1_1} 
-\cdots -2a^2_{n-2}\partial _{ b^1_{n-2}} -a^1_{n-1}\partial _{ b^1_{n-1}}.  
\end{array}  \end{equation}         
\medskip  {\bf Case 2}$\quad n=2k+1$   
\begin{equation} \label{V1o}  \begin{array}{rcl}  
{\bf V}_{1} &=&    
a^1_1\partial_{  a^2_0}+a^1_2\partial_{  a^2_1} +\cdots 
+a^1_{n-1}\partial_{ a^2_{n-2}}   
\\[0.3cm]& & \D + 2a^2_1\partial_{  a^3_0}+2a^2_2\partial_{ a^3_1} 
 +\cdots +2a^2_{n-3}\partial_{ a^3_{n-4}} 
\\ [0.3cm]  & & \dotfill  
\\[0.3cm]& &\D  +\frac{n+1}{2}a^{\frac{n+1}{2}}_1\partial_{ b^{\frac{n+1}{2}}_0}   
\\[0.3cm]& &\D +\frac{n+3}{2}b^{\frac{n-1}{2}}_1 \partial_{ b^{\frac{n-1}{2}}_0} 
 +\frac{n+1}{2}b^{\frac{n-1}{2}}_1\partial_{ b^{\frac{n-1}{2}}_1} 
 +\frac{n-1}{2}a^{\frac{n-1}{2}}_3\partial_{ b^{\frac{n-1}{2}}_2}   
\\[0.3cm] & & \dotfill    
\\[0.3cm] & & \D  +(n-1)b^2_1\partial_{ b^2_0}  +(n-2)b^2_2\partial_{ b^2_1}
+\cdots 2a^2_{n-2}\partial_{ b^2_{n-3}}   \\[0.3cm]&&\D +nb^1_1\partial_{ b^1_0}
 +(n-1)b^1_2\partial_{ b^1_1} +\cdots 
+2b^1_{n-1}\partial_{ b^1_{n-2}} +a^1_n\partial_{ b^1_{n-1}},  
\end{array}  \end{equation}       
\begin{equation} \label{V2o}  \begin{array}{rcl}  
{\bf V}_{2} &=&   
 b^1_1\partial_{ b^2_0}+b^1_2\partial_{ b^2_1} 
+\cdots +b^1_{n-1}\partial_{  a^2_{n-2}}   
\\[0.3cm]  & & \D +2b^2_1\partial_{ b^3_0}+2b^2\partial_{ b^3_1} 
 +\cdots +2b^2_{n-3}\partial _{ a^3_{n-4}}  
\\[0.3cm]  & &  \dotfill  
\\[0.3cm]  & & \D  +\frac{n+1}{2}a^{\frac{n+1}{2}}_1\partial _{ a^{\frac{n+1}{2}}_0}   
\\[0.3cm]  & & \D +\frac{n+3}{2}a^{\frac{n-1}{2}}_1\partial_{ a^{\frac{n-1}{2}}_0} 
 +\frac{n+1}{2}a^{\frac{n+1}{2}}_2\partial_{ b^{\frac{n-1}{2}}_1} 
 +\frac{n-1}{2}a^{\frac{n-1}{2}}_3\partial_{ a^{\frac{n-1}{2}}_2}   
\\[0.3cm]  & & \dotfill     
\\[0.3cm]  & & \D   +(n-1)a^2_1\partial_{ a^2_0} +(n-2)a^2_2\partial_{  a^2_1}
 +\cdots+ 2a^2_{n-2}\partial_{ a^2_{n-3}}   
\\[0.3cm] & & \D  +na^1_1\partial_{ a^1_0} +(n-1)a^1_2\partial_{ a^1_1} 
+\cdots +2a^1_{n-1}\partial_{ a^1_{n-2}} +a^1_n\partial_{ a^1_{n-1}},  
\end{array}  \end{equation}     
\begin{equation} \label{V3o}  \begin{array}{rcl}  
{\bf V}_{3} &=& 
 -na^2_0 \partial_{ a^1_0} -(n-1)a^2_1 \partial_{ a^1_1} 
-\cdots 2a^2_{n-2}\partial_{  a^1_{n-2}} -b^1_{n-1} \partial _{ a^1_{n-1}}   
\\[0.3cm] & & \D -\left[(n-1)a^3_0+a^1_0 \right]\partial_{  a^2_0}   
 -\left[(n-2)a^3_1+a^1_1 \right]\partial_{ a^2_1}   
\\[0.3cm] & & - \cdots -\left[2b^2_{n-3}+a^1_{n-3} \right]\partial_{ a^2_{n-3}}   
\\[0.3cm] & &  \dotfill   
\\[0.3cm] & & \D   -\frac{n-1}{2}[a^{\frac{n-1}{2}}_1 +b^{\frac{n-1}{2}}_1]
 \partial _{ a^{\frac{n+1}{2}}_1} -\cdots -[a^1_{n-2}+b^1_{n-2}]\partial_{ a^2_{n-2}}   
\\[0.3cm] & &\dotfill    
\\[0.3cm] & & \D  -[(n-1)b^3_0+b^1_0 ] \partial _{ b^2_0}  
  - \cdots -[2a^2_{n-3}+b^1_{n-3} ]\partial_{ b^2_{n-3}}   
\\[0.3cm] & &  -nb^2_0 \partial_{ b^1_0} -(n-1)b^2_1 \partial_{  b^1_1}
-\cdots  -2a^2_{n-2}\partial _{ b^1_{n-2}} -a^1_{n-1}\partial_{ b^1_{n-1}}.  
\end{array}  \end{equation}          
 We remark that in both cases the vector fields ${\bf V}_1$, ${\bf V}_2$ and 
${\bf V}_3$  correspond to the generators $\bf T$, $\bf X$ and $\bf H$ given by 
 (\ref{TXH}) respectively. Moreover, it is easy to verify directly that the 
 vector fields $-{\bf V}_1$, $-{\bf V}_2$ and $-{\bf V}_3$ satisfy the 
 same commutator relations (\ref{COMM}) as  $\bf T$, $\bf X$ and $\bf H$.  
We conclude therefore that ${\bf V}_i,$ $i=1,2,3$ represent the infinitesimal 
 action of the isometry group $I(\mathbb{R}_1^2)$ on the parameter space  
$\Sigma$ defined by ${\cal K}^n(\mathbb{R}_1^2)$ for  each $n \ge 1$ and we 
 have proven the ITKT analogue of  Lemma \ref{cayley} of Cayley \cite{AC}:   
\begin{lem} \label{yue} 
 The action of the isometry group $I(\mathbb{R}_1^2)$ on the vector space  
${\cal K}^n(\mathbb{R}_1^2)$ has  the infinitesimal generators (\ref{V1e}),  
(\ref{V2e}) and (\ref{V3e}) when $n$ is even and  (\ref{V1o}), (\ref{V2o})  
and (\ref{V3o}) when $n$ is odd.   
\end{lem}    
\begin{exa} {\rm  Consider the  
vector space ${\cal K}^4(\mathbb R^2_1)$.  
Using the formulas (\ref{K4}), (\ref{V1e}), (\ref{V2e}) and (\ref{V3e}), 
 we derive  the three vector fields representing the infinitesimal action of 
 the isometry group $I(\mathbb{R}_1^2)$ on the vector space 
${\cal K}^4(\mathbb{R}_1^2)$.   
\begin{equation}  \begin{array}{rcl}  
{\bf V}_{1} &=& 
 a^1_1\partial _{a^2_0}   +a^1_2\partial _{a^2_1}+a^1_3\partial _{a^2_2}   
\\ [0.3cm] {}&{}& +2a^2_1\partial _{a^3_0}  
\\ [0.3cm] {}&{}& + 3b^2_1\partial _{b^2_0}+2a^2_2\partial _{b^2_1}  
\\ [0.3cm] {}&{}& + 4b^1_1\partial _{b^1_0}+3b^1_2\partial _{b^1_1} 
 +2b^1_3\partial _{b^1_2}+a^1_4\partial _{b^1_3},  
\end{array}  \end{equation}   
\begin{equation} \begin{array}{rcl}  
{\bf V}_{2} &=& 
 b^1_1\partial _{b^2_0}   +b^1_2\partial _{b^2_1}+b^1_3\partial _{a^2_2}  
 \\ [0.3cm] {}&{}& + 2b^2_1\partial _{a^3_0}  
\\ [0.3cm] {}&{}&+ 3a^2_1\partial _{a^2_0}+2a^2_2\partial _{a^2_1} 
 \\ [0.3cm] {}&{}&+ 4a^1_1\partial _{b^1_0}+3a^1_2\partial _{b^1_1} 
 +2a^1_3\partial _{b^1_2}+a^1_4\partial _{a^1_3},  
\end{array} \end{equation}   
\begin{equation} \begin{array}{rcl} 
 {\bf V}_{3} &=&  
-4a^2_0\partial _{a^1_0}-3a^2_1 \partial_{a^1_1}   
  -2a^2_2\partial _{a^1_2}-b^1_3\partial _{a^1_3}   
\\ [0.3cm] {}&{}& -(3a^3_0+a^1_0)\partial _{a^2_0}
 -(2b^2_1+a^1_1) \partial_{a^2_1}  
\\ [0.3cm] {}&{}& -2(a^2_0+b^2_0)  \partial _{a^3_0}
-(a^1_2+b^1_2)  \partial _{a^2_2}  
\\ [0.3cm] {}&{}& -(3a^3_0+b^1_0)\partial _{b^2_0}-(2a^2_1+b^1_1)\partial _{b^2_1}  
\\ [0.3cm] {}&{}& -4b^2_0\partial _{b^1_0}-3b^2_1\partial _{b^1_1}
  -2a^2_2\partial _{b^1_2}-a^1_3\partial _{b^1_3}.  
\end{array} \end{equation}  
}\end{exa} 
 \bigskip  
 \begin{exa} {\rm  Consider the  vector space ${\cal K}^5(\mathbb R^2_1)$.
 Using the formulas  (\ref{K5}) and  (\ref{V1o}), (\ref{V2o})
 and (\ref{V3o}), we derive  the three vector  fields representing the 
infinitesimal action of the isometry group $I(\mathbb{R}_1^2)$  
on the vector space ${\cal K}^5(\mathbb{R}_1^2)$.   
\begin{equation} \begin{array}{rcl}  
{\bf V}_{1} &=&  
a^1_1\partial _{a^2_0}   +a^1_2\partial _{a^2_1} 
+a^1_3\partial _{a^2_2}+a^1_4 \partial_{a^2_3}   
\\ [0.3cm] {}&{}&+ 2a^2_1\partial _{a^3_0}+2a^2_2 \partial _{a^3_1} 
 \\ [0.3cm] {}&{}&+ 3a^3_1 \partial _{b^3_0}  
\\[0.3cm] {}&{}&+ 4b^2_1\partial _{b^2_0}+3b^2_2\partial _{b^2_1}
 +2a^2_3  \partial_{b^2_2}  
\\ [0.3cm] {}&{}& + 5b^1_1\partial _{b^1_0}+4b^1_2\partial _{b^1_1} 
 +3b^1_3\partial _{b^1_2}+2b^1_4\partial _{b^1_3}+a^1_5  \partial _{b^1_4},  
\end{array} \end{equation}   
\begin{equation} \begin{array}{rcl}  
{\bf V}_{2} &=&  
b^1_1\partial _{b^2_0}   +b^1_2\partial _{b^2_1}
 +b^1_3\partial _{b^2_2}+b^1_4 \partial_{a^2_3}   
\\ [0.3cm] {}&{}&+ 2b^2_1\partial _{b^3_0}+2a^2_2 \partial _{a^3_1}  
\\ [0.3cm] {}&{}& + 3a^3_1 \partial _{a^3_0}  
\\[0.3cm] {}&{}&+ 4a^2_1\partial _{a^2_0}+3a^2_2\partial _{a^2_1}
 +2a^2_3  \partial_{a^2_2}   
\\ [0.3cm] {}&{}&+ 5a^1_1\partial _{a^1_0}+4a^1_2\partial _{a^1_1} 
 +3a^1_3\partial _{a^1_2}+2a^1_4\partial _{a^1_3}+a^1_5  
\partial _{a^1_4},  
\end{array} \end{equation}   
\begin{equation} \begin{array}{rcl}  
{\bf V}_{3} &=&  
-5a^2_0\partial _{a^1_0}   -4a^2_1 \partial _{a^1_1}-3a^2_2 \partial _{a^1_2} 
 -2a^2_3 \partial_{a^1_3}-b^1_4  \partial _{a^1_4}  
\\[0.3cm] {}&{}& -(4a^3_0+a^1_0)  \partial_{a^2_0}-(3a^3_1+a^1_1)  \partial_{a^2_1}
  -(2b^2_2+a^1_2)  \partial_{a^2_2} 
\\ [0.3cm] {}&{}&  -(3b^3_0+2a^2_0)  \partial_{a^3_0}
-2(b^2_1+a^2_1)  \partial_{a^3_1}  -(a^1_3+b^1_3)  \partial_{a^2_3} 
\\ [0.3cm] {}&{}& -(3a^3_0+2b^2_0)  \partial_{b^3_0}  
\\[0.3cm] {}&{}& -(4b^3_0+b^1_0)  \partial_{b^2_0}
-(3a^3_1+b^1_1)  \partial_{b^2_1}  -(2a^2_2+b^1_2)  \partial_{b^2_2}  
\\ [0.3cm] {}&{}& -5b^2_0\partial _{b^1_0}   
-4b^2_1 \partial _{b^1_1} -3b^2_2 \partial _{b^1_2} -2a^2_3 \partial_{b^1_3
}-a^1_4  \partial _{b^1_4}.  
\end{array} \end{equation} 
 }\end{exa}   Using the result of Lemma \ref{yue} we can now employ the 
infinitesimal generators  to compute the $I(\mathbb{R}_1^2)$-invariants.      
\begin{prop}   A function $I:\, \Sigma \rightarrow \mathbb{R}$ is an 
$I(\mathbb{R}_1^2)$-invariant  of the induced action of the isometry group 
$I(\mathbb{R}_1^2)$ on the vector space  ${\cal K}^n(\mathbb{R}_1^2)$ for 
a specific $n\ge 1$ if and only if it satisfies the  infinitesimal criteria  
\begin{equation} \label{crit}  
{\bf V}_1 (I) = {\bf V}_2(I)  = {\bf V}_3 (I) = 0,  
\end{equation}  where $\Sigma$ is the parameter space of  
${\cal K}^n(\mathbb{R}_1^2)$  and the vector fields ${\bf V}_i,$ $i=1,2,3$ 
are given by (\ref{V1e}), (\ref{V2e})  and (\ref{V3e}) when $n$ is even and 
(\ref{V1o}), (\ref{V2o}) and (\ref{V3o})  when $n$ is odd. \label{P2} \end{prop}       
\begin{coro} For a given $n\ge 1$ the parameter $a_n^1$ 
(refer to the formulas (\ref{C2}) and (\ref{C1}) when $n$ is  even and odd respectively) 
is a fundamental $I(\mathbb{R}_1^2)$-invariant of the vector space
 ${\cal K}^n(\mathbb{R}_1^2)$.  \end{coro}   \noindent {\em Proof.} 
Follows from Proposition \ref{P2} and the formulas (\ref{V1e}), (\ref{V2e}) 
 and (\ref{V3e}) when $n$ is even and (\ref{V1o}), (\ref{V2o}) and (\ref{V3o}) 
 when $n$ is odd. \hfill $\Box$    \bigskip   In view of Proposition \ref{P2}, 
the problem  of the determination of the space of $I(\mathbb{R}_1^2)$-invariants
 reduces to  solving the system  of linear PDEs (\ref{crit}). For larger values of
 $n$ the problem  becomes very  challenging computationally. 
The method of characteristics may fail,  in which case one can employ the method of 
undetermined coefficients in conjuncture  with the result of Theorem \ref{FT}, 
as well as  computer algebra.  This technique was  used with a remarkable success in 
 Horwood {\em et al} \cite{HMS}  to solve the problem  of the determination of the 
space of $I(\mathbb{R}^3)$-invariants of the vector space  ${\cal K}^2(\mathbb{R}^3)$,
 where $\mathbb{R}^3$ denotes the Euclidean space.     The concept of a {\em covariant} 
in ITKT was introduced in \cite{SY}.  Proposition \ref{P2} entails 
the corresponding criterion for  $I(\mathbb{R}_1^2)$-covariants of the 
vector spaces ${\cal K}^n(\mathbb{R}_1^2),$  $n \ge 1$.  
\begin{theo} 
Let ${\cal K}^n(\mathbb{R}_1^2)$ be the vector space  of Killing tensors of 
valence $n$ defined in the Minkowski plane $\mathbb{R}_1^2$  
 for a fixed $n \ge 1$. A function 
$C:\, \Sigma\times\mathbb{R}_1^2 \rightarrow \mathbb{R}$ 
 is an $I(\mathbb{R}_1^2)$-covariant of ${\cal K}^n(\mathbb{R}_1^2)$ 
if and only  if it satisfies the infinitesimal invariance conditions  
\begin{equation} \label{crit1}  
\tilde{\bf V}_1 (C) = \tilde{\bf V}_2(C)  = \tilde{\bf V}_3 (C) = 0,  
\end{equation} 
 where the infinitesimal generators are  
\begin{equation} \begin{array}{rcl}  
\tilde{\bf V}_1 & = & {\bf V}_1 + \partial_t, 
\\[0.3cm]  \tilde{\bf V}_2 & = & {\bf V}_2 + \partial_x, 
\\[0.3cm]  \tilde{\bf V}_3 & = & {\bf V}_3 + x\partial_t + t\partial_x, 
\\[0.3cm]  
\end{array} \end{equation}  
$\Sigma$ is the parameter space of  ${\cal K}^n(\mathbb{R}_1^2)$ 
and the vector  fields ${\bf V}_i,$ $i=1,2,3$ are given by (\ref{V1e}), (\ref{V2e}) 
and  (\ref{V3e}) when $n$ is even and (\ref{V1o}), (\ref{V2o}) and (\ref{V3o})
 when  $n$ is odd.  \end{theo}  \section{Conclusions}      
After all, in this paper we have formulated and proven only  {\em an} 
ITKT analogue of Cayley's Lemma in CIT. A similar result for the vector 
 spaces ${\cal K}^n(\mathbb{R}^2)$, $n\ge 1$ (here $\mathbb{R}^2$ denotes the 
 Euclidean plane) can be obtained  {\em mutatis mutandis}. 
Indeed, it is obvious that  the corresponding formulas will differ only by signs.
 More challenging problems are to  extend the result to two-dimensional spaces of
 non-zero curvature, namely when  the  underlying manifold is $\mathbb{S}^2$ 
(two-sphere) or $\mathbb{H}^2$ (hyperbolic plane).  
The work in this direction is underway. 
     
\bigskip    \noindent{\bf Acknowledgements}     
 I am very grateful to my supervisor Professor Roman Smirnov for
 introducing to me ITKT,  formulating the problem, his unfailing support and help 
in preparing the manuscript.  I also thank Professors Dorette Pronk  and Irina Kogan
 for useful discussions  and bringing to our attention  the reference \cite{Sha}
 respectively, as well as  Dr. Alexander Zhalij for bringing to our attention the
 reference \cite{Bo}. Many thanks also to Professor Keith Johnson,  
Joshua Horwood and Dennis The for useful comments and suggestions that have helped 
to improve the presentation of our results. This research was supported by
 an Izaak Walton Killam Memorial Predoctoral Scholarship.                 

\begin{thebibliography}{99}      \bibitem{FO1} M. Fels and P. J. Olver,   Moving coframes.  I. A practical algorithm,  { Acta Appl. Math.} {\bf 51} (1998) 161--213.        \bibitem{FO2} M. Fels and P. J. Olver,    Moving coframes. II. Regularization and theoretical foundations,   {Acta Appl. Math.} {\bf 55} (1999) 127--208.      \bibitem{KO} I. A. Kogan,   Two algorithms for a moving frame construction,   Canad. J. Math. {\bf 55} (2003) 266--291.       \bibitem{Hil} D. Hilbert,   {\it Theory of Algebraic Invariants},   Cambridge University Press, 1993.      \bibitem{PJO} P. J. Olver,   {\em  Classical Invariant Theory},     London Mathematical Society,  Student Texts \textbf{44},    Cambridge University Press, 1999.    \bibitem{MST1} R. G. McLenaghan, R. G. Smirnov and D. The,     Group invariant classification of separable Hamitonian systems    in the Euclidean plane and the $O(4)$-symmetric  Yang-Mills   theories of Yatsun,   J. Math. Phys. {\bf 43} (2002) 1422--1440.             \bibitem{HMS} J. T. Horwood, R. G. McLenaghan and R. G. Smirnov,   Invariant theory and the geometry of orthogonal coordinate webs   in  the Euclidean space, in preparation.      \bibitem{MMS} R. G. McLenaghan, R. Milson and R. G. Smirnov,   Killing tensors as irreducible representations of the general  linear group,   in preparation.       \bibitem{SY} R. G. Smirnov and J. Yue,   Covariants, joint invariants and  the problem of  equivalence   in the invariant theory of Killing tensors,  in preparation.         \bibitem{MST2} R. G. McLenaghan, R. G. Smirnov and D. The,   Towards a classification of cubic integrals of motion,   to  appear in the Proceedings of the conference   ``Superintegrability in Quantum and Classical Mechanics''   (September 16-21, 2002, CRM, Montreal).          \bibitem{DHMS} R. J. Deeley, J. T. Horwood, R. G. McLenaghan and R. G. Smirnov,   Theory of algebraic invariants of vector spaces of Killing tensors:   Methods for  computing the fundamental invariants.   Proceedings of the 5th International Conference   ``Symmetries in Nonlinear Mathematical Physics'',   Part III. Algebras, Groups and Representation Theory, (2004) 1079-1086.       \bibitem{MST3} R. G. McLenaghan, R. G. Smirnov  and D. The,   An extension  of the classical  theory of algebraic invariants   to pseudo-Riemannian geometry and Hamiltonian mechanics,    J. Math. Phys. {\bf 45} (2004) 1079--1120.       \bibitem{MST4} R. G. McLenaghan, R. G. Smirnov  and D. The,   An invariant classification of orthogonal coordinate webs,   Contemp. Math. {\bf 337} (2003) 109--120,Proceedings of the conference   ``Recent advances in Lorentzian and Riemannian geometries''   (January 15-18, 2003, Baltimore).       \bibitem{MST5} R. G. McLenaghan, R. G. Smirnov and D. The,   ``The 1881 problem of Morera revisited,''   {\it Diff. Geom. Appl.,} 2001, 333-241,  in the Proceedings of  ``The 8th Conference on Differential Geometry and Its Applications''   (August 27-31, 2001, Opava, Czech Republic),   Kowalski O., Krupka D. and Slov\'{a}k J eds., 
Silesian University at Opava.        
\bibitem{MST6} R.  G.  McLenaghan, R.  G.  Smirnov and D.  The,   
``Group invariants of Killing tensors in the Minkowski plane'',   
Proceedings of ``Symmetry and Perturbation Theory - SPT2002'',  
 the conference held in Cala Gonone, 19-26 May 2002, 
  S.  Abenda, G.  Gaeta and S.  Walcher eds., World Scientific, 2003, 153--162.        
\bibitem{GCR} G.-C. Rota,   
{\em Indiscrete Thoughts}, Brikh\"{a}user, 1997.     
\bibitem{BO} M. Boutin,   
{\em On Invariants of Lie Group Actions and their Application to  
some Equivalence Problems},   
PhD thesis, University of  Minnesota, 2001.       
\bibitem{BMS} A. T. Bruce, R. G. McLenaghan and R. G. Smirnov,  
 A geometrical approach to the problem of integrability of  Hamiltonian systems 
 by separation of variables,   
J. Geom. Phys. {\bf 39} (2001) 301--322.     
\bibitem{PJO1} P. J. Olver,    
{\em Equivalence, Invariants and Symmetry},   Cambridge University Press, 1995.     
\bibitem{AC} A. Cayley,   
A second memoir on quantics,     
 Phil. Trans. Roy. Soc. London, {\bf 144} (1856) 561--578.      
\bibitem{AT} A. C. Thompson, 
{\it Minkowski Geometry}, Encyclopedia of Mathematics and its Applications {\bf 63}, 
 Cambridge University Press, Cambridge, 1996. 
\bibitem{Ri} B. Riemann,   
{\em \"{U}ber die Hypothesen welche der Geometrie zu Grunde liegen},  
 the lecture given on 10 June 1854 at G\"{o}ttingen University.        
 \bibitem{FK72} F. Klein,    
Vergeleichende Betrachtungen \"{u}ber neuere geometrische Forshcungen.  
 Erlangen: A. Duchert (1872).    
\bibitem{FK93} F. Klein,   
Vergeleichende Betrachtungen \"{u}ber neuere geometrische Forshcungen.  
 Math. Ann. {\bf 43} (1893) 63--100 (Revised version of \cite{FK72}).     
\bibitem{Car} \'{E}. Cartan,  
 {\em La M\'{e}thode du Rep\`{e}re Mobile, la Th\'{e}orie 
   des Groups Continues, et les Espaces G\'{e}n\'{e}ralis\'{e}s.}   
Expos\'{e}s de G\'{e}om\'{e}trie {\bf 5},  Hermann, Paris, 1935.      
\bibitem{Sha} R. W. Sharpe, {\em Differential Geometry:  
Cartan's Generalization of Klein's Erlangen Program},  Springer, 1997.  
\bibitem{Arv} A. Arvanitoyeorgos, 
{\em An Introduction to Lie Groups and the Geometry of Homogeneous Spaces}, 
 Student Mathematical Library {\bf 22}, AMS, 2003.       
 \bibitem{Sch}J. A. Schouten,   
 \"{U}ber Differentalkomitanten zweier kontravarianter Gr\"{o}ssen,  
 { Proc. Kon. Ned. Akad. Amsterdam} {\bf 43} (1940) 449--452.     
\bibitem{De} R. P. Delong, Jr.,  
 {\it Killing tensors and the Hamilton-Jacobi equation},   
  PhD thesis, University of Minnesota, 1982.    
 \bibitem{Do} P. Dolan, A. Kladouchou, C. Card,  
 On the significance of Killing tensors,  
 {Gen. Relativity and Gravitation}, {\bf 21} (1989) 427--437.      
\bibitem{Be} S. Benenti,  
 Separability in Riemannian manifolds,  
to appear in Proc. Roy. Soc.       
\bibitem{BM} A. V. Bolsinov and V. S. Matveev,   
Geometrical interpretation of Benenti systems,    
J. Geom. Phys. {\bf 44} (2003) 489--506.      
\bibitem{Cr} M. Crampin,   Conformal Killing tensors with vanishing torsion
 and the separation   of variables in the Hamilton-Jacobi equation,  
 Diff. Geom. Appl. {\bf 18} (2003) 87--102.      
\bibitem{Ei1} L. P. Eisenhart,   
Separable systems of St\"{a}ckel,   
Ann. Math. {\bf 35} (1934) 284--305.        
\bibitem{Ei2} L. P. Eisenhart,  
 St\"{a}ckel systems in conformal Euclidean space,   
Ann. Math. {\bf 36} (1934) 57--70.      
\bibitem{FN} W. I. Fushchich. and A. G. Nikitin,  
 {\em Symmetries of Equations of Quantum Mechanics},   
Allerton Press Inc., New York, 1994.     
\bibitem{Kal} E. Kalnins,   
On the separation of variables for the Laplace equation   
in two- and three-dimensional Minkowski space,   
{SIAM J. Math. Anal.} {\bf 6} (1975) 340--374.     
\bibitem{Ka} E. G. Kalnins,   
{\it Separation of Variables for Riemannian Spaces of Constant Curvature}  
 Longman Scientific $\&$ Technical, New York, 1986.       
\bibitem{KM1} E. G. Kalnins and W. Miller, Jr.,  
 Killing tensors and nonorthogonal variable separation for Hamilton-Jacobi equations,   
 SIAM J. Math. Anal. {\bf 12} (1981) 617--629.       
\bibitem{KM2} E. G. Kalnins and W. Miller, Jr.,   
Conformal Killing tensors and variable separation for Hamilton-Jacobi equations,   
SIAM J. Math.  Anal. {\bf 14} (1983) 126--137.      
\bibitem{Mi} W. Miller, Jr.,   
{\it Symmetry and Separation of Variables}  
 Addison-Wesley, New York, 1977.      
\bibitem{MF} O. I. Mokhov and E. V. Ferapontov,  
 Hamiltonian pairs generated by skew-symmetric Killing tensors  
 on spaces of constant curvature,   
Funct. Anal. Appl. {\bf 28} (1994) 123--125.       
\bibitem{Ta} M. Takeuchi, 
 Killing tensor fields on spaces of constant curvature,   
{ Tsukuba J. Math.} {\bf 7} (1983) 233--255.     
\bibitem{Th} G. Thompson,   
Killing tensors in spaces of constant curvature,   
{ J. Math. Phys.} {\bf 27} (1986) 2693--2699.     
\bibitem{Bo} M. B\^ocher,   
{\em \"{U}ber die Reihenentwickelungen der Potentialtheorie}  
 (mit einem Vorwort von Felix Klein), Leipzig, Teubner Verlag, 1894.      
\bibitem{Dar} G. Darboux,  
 {\em Le\c{c}ons sur les Syst\`{e}mes Orthogonaux et les Coordonn\'{e}es Curvilingnes}, 
  Paris, Gauthier-Villars, 1910.     
\bibitem{Sm} R. G. Smirnov,  
 Bi-Hamiltonian formalism: A constructive approach,  
 Lett. Math. Phys. {\bf 41} (1997) 333--347.      
\bibitem{Ser} A. Sergyeyev,   
A simple way of making a Hamiltonian system into a bi-Hamiltonian one,  
 to appear in Acta Appl.  Math.      
\end{thebibliography}
\end{document}